\documentclass[a4paper, 11pt]{amsart}

\usepackage[latin1]{inputenc}
\usepackage{amssymb, amsmath}
\usepackage{mathrsfs}
\usepackage{amsthm}
\usepackage[english]{babel}

\usepackage[numeric, lite, initials, nobysame]{amsrefs}
\input xy
\xyoption{arrow} \xyoption{matrix} \xyoption{frame}%\OnlyOutlines
\xyoption{curve}

\let\oldmarginpar\marginpar
\renewcommand\marginpar[1]{\-\oldmarginpar[\raggedleft\footnotesize #1]%
{\raggedright\footnotesize #1}}

\def\:{\colon}
\def\<{\left\langle}
\def\>{\right\rangle}
\def\({\left(}
\def\){\right)}
\def\epsilon{\varepsilon}
\def\phi{\varphi}
\def\subset{\subseteq}

\def\leq{\leqslant}
\def\geq{\geqslant}

\def\lra{\longrightarrow}
\def\Lra{\Longrightarrow}
\def\mapsto{\longmapsto}
\def\smash{\wedge}
\def\xra{\xrightarrow}

\def\op{\text{op}}

\def\incl{\text{incl}}

\def\wh{\widehat}
\def\odd{\text{odd}}

\def\DR{\mathscr{D}_R}
\def\Cl{\mathcal C \ell}
\def\SS{\mathbb S}\def\C{\mathbb C}
\def\H{{H}}
\def\Z{\mathbb Z}
\def\F{\mathbb F}

\newtheorem{thm}{Theorem}
\newtheorem{lem}[thm]{Lemma}
\newtheorem{prop}[thm]{Proposition}
\newtheorem{cor}[thm]{Corollary}

\newtheorem*{thm*}{Theorem}
\newtheorem*{cor*}{Corollary}
\newtheorem*{prop*}{Proposition}

\numberwithin{equation}{section}
\numberwithin{thm}{section}

\theoremstyle{remark}
\newtheorem{rem}[thm]{Remark}

\theoremstyle{definition}
\newtheorem{defn}[thm]{Definition}
\newtheorem{exmp}[thm]{Example}

\DeclareMathOperator{\Hom}{Hom}

\DeclareMathOperator{\hocolim}{hocolim} 
\DeclareMathOperator{\colim}{colim} 
 
\DeclareMathOperator{\T}{T} 
\DeclareMathOperator{\DDer}{\mathscr Der} 
\DeclareMathOperator{\Tor}{Tor}
\DeclareMathOperator{\Der}{Der}

\newcommand{\ie}{i.e.}

\newcommand{\eg}{e.g.}

\title{Clifford algebras from quotient ring spectra}

\author{A. Jeanneret}
\address{Mathematisches Institut, Sidlerstrasse 5, 3012 BERN, Switzerland}
\email{alain.jeanneret@math.unibe.ch}

\author{S. W\"uthrich}
\address{SBB, Br\"uckfeldstrasse 16, 3000 BERN, Switzerland}
\email{samuel.wuethrich@sbb.ch}

\date{18.\@ January 2011}

\begin{document}

\begin{abstract}
We give natural descriptions of the homology and cohomology algebras of regular quotient
ring spectra of even $E_\infty$-ring spectra. We show that the homology is a Clifford
algebra with respect to a certain bilinear form naturally associated to the quotient ring
spectrum $F$. To identify the cohomology algebra, we first determine the derivations of
$F$ and then prove that the cohomology is isomorphic to the exterior algebra on the
module of derivations. We treat the example of the Morava $K$-theories in detail.
\end{abstract}

\maketitle

\section{Introduction}

It has long been a difficult problem to realize quotient constructions in stable homotopy
theory. The situation changed completely with the introduction of point-set categories of
spectra endowed with monoidal structures, for instance in \cite{ekmm}. Since then, the
definition of a large class of quotient constructions has become a pure formality.
Namely, suppose that $R$ is an $E_\infty$-ring spectrum and that $I\subset \pi_*(R)= R_*$
is an ideal of the homotopy ring of $R$ generated by a regular sequence. Then there is a
spectrum $F$ equipped with a map $R\to F$ which induces an isomorphism $F_*\cong R_*/I$.
Moreover, $F$ is unique up to equivalence, see Remark \ref{htyregquot}.

Such regular quotients of $R$ arise naturally as objects in the derived category of
$R$-module spectra $\mathscr D_R$. Working in this category makes it much easier to study
multiplicative structures. Partly, this is due to the fact that $\DR$ is equipped with a
monoidal structure, induced by the smash product $\smash_R$. In particular, Strickland
\cite{strickland} showed that
a regular quotient can always be realized as an $R$-ring
spectrum, \ie\ as a monoid in $\DR$, if $R_*$ forms a domain and $R$ is even, meaning
that $R_*$ is trivial in odd degrees.

A fundamental problem is to compute the homology and cohomology
algebras of regular quotients $F_*^R(F)=\pi_*(F\smash_R F)$ and
$F^*_R(F)=\DR^*(F,F)$, respectively. Whereas the underlying graded
$F_*$-modules are trivial to determine if $R$ is even, the
multiplicative structures have only been identified in special
cases up to now,
see \cite{angeltviet}, \cite{b-j}, \cite{b-l}, \cite{laz}
and \cite{strickland}. The main goal of this article is to determine the
homology and cohomology algebras in general. Our descriptions are
valid for arbitary products on $F$ and functorial in nature. In
particular, they are independent on any choices, such as the
fixing of generators of $I$. This is important in \cite{jw}, where
the results proved here  are used to solve the classification
problem of
$R$-ring structures on regular quotients.

%Our methods do not seem to depend on any particular features of
%the homotopy category of module spectra $\DR$ of an
%$E_\infty$-ring $R$. We expect that they can be carried over to
%many other settings, \eg\ the derived categories of ordinary
%commutative rings. We hope to return to this in future work.

\smallskip

We do not restrict to regular quotient rings, but consider arbitrary quotient rings of an
even $E_\infty$-ring $R$, \ie\ $R$-rings $F$ with $F_*\cong R_*/I$ for some ideal
$I\subset R_*$. We write $F=R/I$ for such an $F$. We study the homology and cohomology of
$F$ with respect to any quotient $R$-ring spectrum $k$ which comes with a unital map
$\pi\: F\to k$. We call $(F, k, \pi)$ with these properties an admissible pair. An
important example of an admissible pair is given by $(F, k, \pi)$, where $F$ is a
quotient ring, $k=F$ as an $R$-module, but endowed with a possibly different product, and
where $\pi$ is the identity map $1_F$.

Our arguments are based on a canonical homomorphism of $k_*$-modules, the {\em
characteristic homomorphism}
\[
\phi\:  k_*\otimes_{F_*} I/I^2[1] \lra k_*^R(F).
\]
Here, $I/I^2[1]$ denotes the graded $F_*$-module $I/I^2$ with degrees raised by one. We
show that $\phi$ is independent of the products on $F$ and $k$ and functorial in both $F$
and $k$. We then use $\phi$ to define the {\em characteristic bilinear form}
\[
b\: \bigl(k_*\otimes_{F_*} I/I^2[1] \bigr) \otimes_{k_*} \bigl(k_*\otimes_{F_*}I/I^2[1]
\bigr) \lra k_*,
\]
Letting $q\: k_*\otimes_{F_*} I/I^2[1]\to k_*$ be the associated quadratic form and
writing $\Cl(k_*\otimes_{F_*} I/I^2[1], q)$ for the Clifford algebra with respect to $q$,
we prove:

\begin{thm*}
For an admissible pair $(F=R/I,k,\pi)$, the characteristic homomorphism
%\[
%\phi\:  k_*\otimes_{F_*} I/I^2[1] \lra k_*^R(F)
%\]
lifts to a natural homomorphism of $k_*$-algebras
\[
\Phi\:  \Cl(k_*\otimes_{F_*} I/I^2[1], q) \lra k^R_*(F).
\]
If $F$ is a regular quotient, then $\Phi$ is an isomorphism.
\end{thm*}

We show that the characteristic bilinear form of $(F, F^\op, 1_F)$ is trivial, where
$F^\op$ denotes the opposite ring of $F$. This leads to a new proof of the fact that
$F_*^R(F^\op)\cong \Lambda(I/I^2[1]))$ is an exterior algebra \cite{laz}. If $F$ is a
{\em diagonal} regular quotient, \ie\ the smash product of quotient rings of the form
$R/x_i$ with $x_i\in R_*$, the characteristic bilinear form of $F$ is diagonal. We prove
that the diagonal elements are determined by the commutativity obstructions of $R/x_i$
introduced in \cite{strickland}.

In the second part of the article, we consider the cohomology modules $k^*_R(F)$, endowed
with the profinite topology. We first show:

\begin{prop*}
For an admissible pair $(F=R/I, k, \pi)$, there exists a natural continuous homomorphism
of $k_*$-modules
\[
\Psi\: k^*_R(F) \lra \Hom_{F_*}^*(\Lambda(I/I^2[1]), k_*).
\]
If $F$ is a regular quotient ring, $\Psi$ is a homeomorphism.
\end{prop*}

For the determination of the cohomology algebra of a regular quotient ring $F$, we
consider the group of derivations $\DDer_R^*(F,F)$. More generally, we study
$\DDer_R^*(F,k)$ for any multiplicative admissible pair $(F,k,\pi)$, \ie\ one for which
$\pi$ is multiplicative. It inherits from $k^*_R(F)$ a linear topology.

\begin{prop*}
For a multiplicative admissible pair $(F=R/I, k, \pi)$ such that $F$ and $k$ are regular
quotients, there is a natural homeo\-morphism
\[
\psi\: \DDer_R^*(F, k) \lra \Hom^*_{F_*}(I/I^2[1], k_*).
\]
\end{prop*}

We then describe the cohomology algebra of $F$ in terms of its derivations. Let
$\wh\Lambda(\DDer^*_R(F, F))$ denote the completed exterior algebra on $\DDer^*_R(F, F)$.

\begin{thm*}
For a regular quotient ring $F=R/I$, there is a canonical homeo\-morphism of
$F^*$-algebras $F^*_R(F) \cong \wh\Lambda(\DDer^*_R(F, F))$.
\end{thm*}

These two statements are generalizations of results of Strickland \cite{strickland}. He
considered the special case where $F$ is a diagonal regular quotient ring of $R$.

In the last section, we discuss the case of the Morava
$K$-theories $K(n)$. We determine explicitly the bilinear form
$b_{K(n)}$. The reader will find in \cite{jw} more examples of
computations of characteristic bilinear forms.

\smallskip

Here is an overview over the contents of this article. In Section \ref{bilinearform}, we
recall some background material from \cite{strickland}, construct the characteristic
homomorphism and characteristic bilinear form of admissible pairs and compute them in
special cases. In Section \ref{homology}, we consider the homology of admissible pairs.
In Section \ref{cohomology}, we study derivations and the cohomology of admissible pairs.
Finally, in Section \ref{examples}, we discuss the example of the Morava $K$-theories.

\subsection{Notation and conventions}\label{notations}
For definiteness, we work in the framework of $\SS$-modules of
\cite{ekmm}. In this setting, $E_\infty$-ring spectra correspond
to commutative $\SS$-algebras. Throughout the paper, $R$ denotes
an even commutative $\SS$-algebra. We also assume that the
coefficient ring $R_*$ of $R$ is a domain (see Remark
\ref{whydomain} for an explanation). Associated to $R$ is the
homotopy category $\DR$ of $R$-module spectra. For simplicity, we
refer to its objects as $R$-modules. The smash product $\smash_R$
endows $\DR$ with a symmetric monoidal structure. We will
abbreviate $\smash_R$ by $\smash$ throughout the paper.

Monoids in $\DR$ are called $R$-ring spectra or just $R$-rings. Unless otherwise
specified, we use the generic notation $\eta_F\: R\to F$ and $\mu_F\: F\smash F\to F$ for
the unit and the multiplication maps of an $R$-ring $F$. Mostly, $\eta_F$ will be clear
from the context, in which case we call a map $\mu_F\: F\smash F\to F$ which gives $F$
the structure of an $R$-ring an $R$-product or just a product. We denote the opposite of
an $R$-ring $F$ by $F^\op$. Its product is given by $\mu_{F^\op}= \mu_F\circ\tau$, where
$\tau \: F\smash F\to F\smash F$ is the switch map. An $R$-ring $(F, \mu_F,\eta_F)$
determines multiplicative homology and cohomology theories $F_*^R(-)=\pi_*(F\smash -) =
\DR^{-*}(R,F\smash -)$ and $F^*_R(-)=\DR^*(-, F)$, respectively, on $\DR$.

Since we are working with non-commutative $R$-rings, we must
carefully describe the various module structures involved. For an
$R$-ring $k$ and an $R$-module $M$, the homology $k_*^R(M)$ is a
$k_*$-bimodule in a natural way. Even if $k_*$ is commutative, the
left and right $k_*$-actions may well be different. However, if we
assume that $k$ is a quotient of $R$, by which we mean that the
unit map $\eta_k\: R\to k$ induces a surjection on homotopy groups
(see Definition \ref{defquot} below), the left and right
$k_*$-actions agree. In this case, we can refer to $k_*^R(M)$ as a
$k_*$-module without any ambiguity. A similar discussion applies
to cohomology $k_R^*(M)$.

We will assume that $k$ is a quotient of $R$ for the rest of this
section.

For $R$-modules $M$ and $N$, we write
\[
\kappa_k \: k_*^R(M) \otimes_{k_*} k_*^R(N) \lra k_*^R(M\smash N)
\]
for the K\"unneth homomorphism, a homomorphism of $k_*$-modules.
Note that $k$ is not required to be commutative for the definition
of $\kappa_k$ (see \cite{sw}*{\S 2}). If $k_*^R(M)$ or $k_*^R(N)$
is $k_*$-flat, then $\kappa_k$ is an isomorphism of $k_*$-modules.

Let $F$ be a second $R$-ring. The composition
\[
m^{k}_{F} \: k_*^R(F) \otimes_{k_*} k_*^R(F) \xra{\kappa_{k}}
k_*^R(F\smash F) \xra{k_*^R(\mu_F)} k_*^R(F)
\]
defines a (central) $k_*$-algebra structure on $k_*^R(F)$, where the unit is given by
$(1_k\smash \eta_F)_* \: k_* \to k_*^R(F)$. In unambiguous situations, we will write
$a\cdot b$ for $m^{k}_{F}(a\otimes b)$.

To relate the homology $k_*^R(M)$ and cohomology $k_R^*(M)$, we
will use the Kronecker duality morphism
\[
d \: k_R^*(M) \lra \Hom^*_{k_*}(k^R_*(M), k_*),
\]
which associates to $f\: M\to k$ the homomorphism of $k^*$-modules
$d(f)=(\mu_k)_* k^R_*(f)$. If $k_*^R(M)$ is $k_*$-free, $d$ is an
isomorphism. This implies that the Hurewicz homomorphism
\[
k^R_*(-)\:  k_R^*(M) \lra \Hom^*_{k_*}(k_*^R(M), k_*^R(k))
\]
is injective whenever $k_*^R(M)$ is $k_*$-free. See \eg\ \cite{sw}*{Lemma 6.2} for a
detailed discussion, which covers in particular the case where $k$ is non-commutative.

For $R$-modules $M$ and $N$, we write $\zeta\: M_* \otimes_{R_*} N_* \lra (M\smash N)_* $
for the canonical map, which is natural in the following graded sense. Two maps of
$R$-modules $f\: \Sigma^k M \to M'$ and $g\: \Sigma^l N \to N'$ induce commutative
diagrams
\begin{equation}\label{gradedcommutative}
\begin{array}{c}
\xymatrix{ M_m \otimes N_n \ar[r]^-{\zeta} \ar[d]_-{f_*\otimes g_*} & (M\smash N)_{m+n}
\ar[d]^-{(-1)^{m\cdot l} f\smash g}
\\
M'_{k+m} \otimes N'_{l+n} \ar[r]^-{\zeta} & (M'\smash N')_{k+m+l+n}.}
\end{array}
\end{equation}

We write $M_*[d]$ for the $d$-fold suspension of a graded abelian group $M_*$, so
$(M_*[d])_k = M_{k-d}$. With this convention, we have $(\Sigma^d M)_* = M_*[d]$ for an
$R$-module $M$. We denote the image of some element $\alpha\in M_k$ under the shift
isomorphism $M_*\cong M_*[d]$ by $\alpha[d]\in (M_*[d])_{k+d}$. We use the convention
$M^*=M_{-*}$. If the ground ring is clear from the context, we omit it from the tensor
product symbol $\otimes$ from now on.

\subsection*{Acknowledgments}
The first author would like to thank Andrew Baker. The present article generalizes
results announced in \cite{b-j}. The second author has benefitted from the stimulating
atmosphere at the IGAT at the EPF Lausanne during the work on this article. He would like
to thank Prof.\@ Kathryn Hess for her support. Furthermore, he would like to express his
gratitude to the Mathematical Institute at the University of Berne for offering him the
opportunity of an inspiring research visit in summer 2006.

\section{The characteristic bilinear form}\label{bilinearform}

\subsection{Quotient modules, quotient rings}\label{regularquotient}

The point of this subsection is to introduce some convenient terminology and to recall
some basic constructions in the category $\DR$, mainly from \cite{strickland}.

\begin{defn}\label{defquot}
A {\em quotient module}\/ of $R$ is an $R$-module $F$ with a map
of $R$-modules $\eta_F\: R\to F$ which induces a surjection on
homotopy groups, that is $F_* \cong R_*/I$ where $I\subset R_*$ is
an ideal. A morphism $f\: F\to G$ of quotient modules of $R$ is a
map of $R$-modules such that $f\circ \eta_F = \eta_G$.
\end{defn}

Let $F$ be a quotient module of $R$ with $F_* =R_*/I$  and let $X$
be the homotopy fibre of $\eta_F\: R\to F$. As the canonical map
$X\to R$ induces an isomorphism $X_*\cong I\subset R_*$, we write
$I$ for $X$. So we have a cofibre sequence of the form
\begin{equation}\label{defineI}
I \xra{\iota} R \xra{\eta_F} F \xra{\beta} \Sigma I.
\end{equation}
We will write $F=R/I$ in the sequel.

Recall that a graded $R_*$-module $F_*$ is said to be a (finite)
regular quotient of $R_*$ if it is isomorphic to $R_*/(x_1, x_2,
\ldots)$ for some (finite) regular sequence $(x_1, x_2, \ldots)$
in $R_*$. There is the following analogous topological notion.

\begin{defn}
A quotient module $F=R/I$ of $R$ is a (finite) {\em regular
quotient module}\/ of $R$ if the ideal $I$ is generated by some
(finite) regular sequence $(x_1, x_2, \ldots)$ in $R_*$.
\end{defn}

We now recall the definition of the building blocks of regular
quotients of $R$. The coefficient ring $R_*$ may be canonically
identified with the graded endomorphisms of $R$ in $\DR$. If $x$
is a given element of $R_d$, we write $R/x$ for the homotopy
cofibre of $x \: \Sigma^d R \to R$. As $R_*$ lies in even degrees,
$R/x$ is well defined up to canonical homotopy equivalence. By
construction, $R/x$ fits into a cofibre sequence of the form
\begin{equation}\label{cofx}
\Sigma^d R \xra{x} R \xra{\eta_x} R/x \xra{\beta_x} \Sigma^{d+1}
R.
\end{equation}
Since $R_*$ is a domain, $(R/x)_* \cong R_*/(x)$.

\begin{rem}\label{htyregquot}
If $F=R/I$ is a regular quotient and $(x_1, x_2, \ldots)$ is a regular sequence
generating $I$, then $F$ is isomorphic in $\DR$ to
\[
R/x_1\smash R/x_2 \smash \cdots := \hocolim_k R/x_1\smash \cdots
\smash R/x_k.
\]
Due to the lack
of a specific reference,
we give a brief
outline of the argument underlying
the proof. We construct by
induction, using
\cite{ekmm}*{V.1, Lemma 1.5}
factorizations
$R/x_1\smash \cdots  \smash R/x_k \to F$
of the unit $\eta\: R \to F$
and from these a map
$\bar\eta\: \hocolim_k R/x_1\smash \cdots
\smash R/x_k \to F$.
By construction,
$\bar\eta$ induces an isomorphism
on homotopy groups and is
thus an isomorphism in $\DR$.
\end{rem}

\begin{defn}
A (regular) {\em quotient ring}\/ of $R$ is an $R$-ring $(F,
\mu_F, \eta_F)$ such that $(F, \eta_F)$ is a (regular) quotient
module of $R$.
\end{defn}

Products on regular quotients of the form $R/x$ have been studied
in \cite{strickland}*{Section 3}.

\begin{prop}\label{basicproducts}
Let $x \in R_d$. If $u$ is in $R_{2d+2}/x$ and $\mu$ is a product
on $R/x$, then $\mu+u\circ(\beta_x\smash\beta_x)$ is another
product. This construction gives a free transitive action of the
group $R_{2d+2}/x$ on the set of products on $R/x$.
\end{prop}

\begin{prop}\label{obstruction}
There is a natural map $c$ from the set of products on $R/x$ to
$R_{2d+2}/x$ such that $c(\mu)=0$ if and only if $\mu$ is
commutative.
\end{prop}

Recall that $R$-ring maps $f\: A\to C$ and $g\: B\to C$ are said
to commute if $ \mu_C \circ (f\smash g) = \mu_C \circ \tau
\circ(f\smash g)\: A\smash B\to C$.

\begin{rem}\label{prodstrickland}
Let $F=R/I$ be a regular quotient of $R$ and $(x_1,x_2,\ldots)$ a
regular sequence generating the ideal $I$. For any products
$\mu_i$ on $R/x_i$, $i\geq 1$, \cite{strickland}*{Prop.\@ 4.8}
implies that there is a unique product on $F=R/I$ such that the
natural maps $R/x_i \to F$ are multiplicative and commute. See
Proposition \ref{diagonalprop} for a generalization.
\end{rem}

\begin{defn}\label{diagonal1}
We call $F$, endowed with the product described in Remark
\ref{prodstrickland}, the {\em smash ring spectrum}\/ of the
$R/x_i$. If we need to be more precise, we refer to the product
map $\mu_F$ on $F$ as the {\em smash ring product}\/ of the
$\mu_i$.
\end{defn}

For the next definition, recall that two $R$-ring spectra $F$ and $G$ are called
equivalent if there is an isomorphism $f\: F\to G$ in $\DR$ which is multiplicative.

\begin{defn}
We call a regular quotient ring $F$ of $R$ {\em diagonal}\/ if it is the smash ring
spectrum of ring spectra $R/x_i$, where $(x_1,x_2,\ldots)$ is a regular sequence. We say
that $F$ is {\em diagonalizable}\/ if it is equivalent to a diagonal regular quotient
ring.
\end{defn}

\begin{cor}\label{realizereg}
Any regular quotient ring of $R_*$ can be realized as the
coefficient ring of a diagonal $R$-ring.
\end{cor}

\begin{proof}
Let $F_*=R_*/(x_1, x_2, \ldots)$ be a regular quotient of $R_*$. By Remark
\ref{htyregquot}, the $R$-module $F=R/x_1\smash R/x_2\smash\cdots$ satisfies
$\pi_*(F)=F_*$. By Proposition \ref{basicproducts}, every $R/x_i$ admits a product.
Finally, endow $F$ with the induced smash ring product.
\end{proof}

\begin{rem}\label{whydomain}
Note that the proof requires each of the elements $x_k$ of the
regular sequence to be a non-zero divisor. This is guaranteed by
our assumption that $R_*$ is a domain.
\end{rem}

Let $(R/x, \mu, \eta)$ be a regular quotient ring, $x\in R_d$, and
$A$ an even $R$-ring. Clearly, there is a unital map $j\: R/x\to
A$ if and only if $x$ maps to zero in $A_*$, and $j$ is unique if
it exists. We will extensively use the following fact:

\begin{prop}\label{natargument}
Let $A$ and $x$ be as above and assume that $A$ is a quotient ring of $R$. Then there
exists a product on $R/x$ such that the canonical map $j\: R/x \to A$ is multiplicative.
\end{prop}

\begin{proof}
Choose an arbitrary product $\mu$ on $F=R/x$. By Proposition
\ref{basicproducts}, any other product on $F$ is of the form
$u\cdot\mu:=\mu+u(\beta \smash \beta)$, for some $u\in
R_{2d+2}/x$. By Proposition \cite{strickland}*{Prop.\@ 3.15},
there is an obstruction $d_F(u\cdot\mu)\in A_{2d+2}$ which
vanishes if and only if $j$ is multiplicative for $u\cdot\mu$.
Furthermore, $d_F(u\cdot\mu)$ is related to $d_F(\mu)$ by
$d_F(u\cdot\mu)= d_F(\mu)+j_*(u)$. Thus, on choosing $u$ with
$j_*(u)=-d_F(\mu)$, $j$ is multiplicative with respect to
$u\cdot\mu$.
\end{proof}

\begin{cor}\label{commdiagonalize}
Let $F=R/I$ be a commutative regular quotient ring of $R$. Then $F$ is diagonalizable.
\end{cor}

\begin{proof}
Let $(x_1, x_2, \ldots)$ be a regular sequence which generates
$I$. By Proposition \ref{natargument}, there are products $\mu_i$
on $R/x_i$ such that the canonical maps $j_i\: R/x_i \to F$ are
multiplicative. By commutativity of $F$, the $j_i$ commute with
each other. By \cite{strickland}*{Prop.\@ 4.8}, they therefore
induce a multiplicative equivalence $j\: \bigwedge_{i\geq 1} R/x_i
\to F$.
\end{proof}

\subsection{Admissible pairs}\label{ac}
In this subsection we introduce the category of admissible pairs.
It will play a central role in the sequel.
\begin{defn}
An {\em admissible pair} is a triple $(F, k, \pi)$ consisting of
two quotient $R$-rings $(F, \mu_F, \eta_F)$,  $(k, \mu_k, \eta_k)$
and a unital map of $R$-modules $\pi\: F\to k$, \ie\ an
$R$-morphism such that $\pi \circ \eta_F = \eta_k$. If $\pi$ is a
map of $R$-ring spectra, we call $(F, k, \pi)$ a {\em
multiplicative admissible pair}.
\end{defn}

Note that $\pi_*\: F_* \to k_*$ is always a ring homomorphism,
even for non-multiplicative admissible pairs, as $(\eta_F)_*\:
R_*\to F_*$ is surjective. We may therefore view $k_*$ as an
$F_*$-module.

\begin{rem}
If $F=R/I$ and $k$ are quotient $R$-rings, a necessary condition
for the existence of a map $\pi$ making $(F,k,\pi)$ into an
admissible pair is that $(\eta_k)_*(I)=0$. If $F$ is a {\em
regular} quotient ring, this condition is sufficient, by
\cite{strickland}*{Lemma 4.7}. If $F=R/x$, the map $\pi$ is
unique.
\end{rem}

Admissible pairs $(F,k, \pi)$ form the objects of a category. The
morphisms between two admissible pairs $(F=R/I, k, \pi)$ and
$(G=R/J, l, \pi')$ are pairs of $R$-ring maps $(f\: F\to G,\; g\:
k\to l)$ which make the diagram
\begin{equation}\label{squarediagram}
\begin{array}{c}
\xymatrix{  F \ar[r]^-{\pi} \ar[d]^-f & k \ar[d]^-{g}
\\
G \ar[r]^-{\pi'} & l}
\end{array}
\end{equation}
commutative. Observe that in this case $I \subset J$. If we say that a certain
construction is ``natural in $F$ and $k$'', we mean that it is a functor on this
category. Similarly, we refer to a morphism as being ``natural in $F$ and $k$'' if it
defines a natural transformation of functors defined on this category.

\begin{exmp}\label{admcoupleFF}
An important example of an admissible pair $(F,k,\pi)$ is the special case where $k$
coincides with $F$ as an $R$-module and $\pi=1_F$, the identity on $F$, but where we
distinguish two products $\mu$ and $\nu$ on $F$.
\end{exmp}

\subsection{The characteristic homomorphism}\label{phi}

Let $(F=R/I, k, \pi)$ be an admissible pair. We define a
homomorphism of $F_*$-modules
\begin{equation}\label{natmap}
\phi^k_F\:  I/I^2[1] \lra k_*^R(F),
\end{equation}
which is natural in $F$ and $k$. Here, we view $k_*^R(F)$ as an
$F_*$-module via the ring homomorphism $\pi_* \: F_* \to k_*$ and
the $k_*$-module structure on $k_*^R(F)$ as discussed in Section
\ref{notations}.

Applying $k\smash -$ to the cofibre sequence \eqref{defineI}
yields
\begin{equation}\label{defsigma}
k\xra{1\smash\eta_F} k\smash F\xra{1\smash\beta} k\smash\Sigma I
\xra{1\smash\iota} \Sigma k.
\end{equation}
We consider the map $\psi\: k\smash F\to k$, defined as
$\psi=\mu_k \circ (1 \smash \pi)$, and we note that $\psi$ is a
retraction of $1\smash\eta_F$, natural in $F$ and $k$. Observe
that if $k=F$ as $R$-modules (Example \ref{admcoupleFF}), then
$\psi$ is just the second product $\nu$ of $F$.

The cofibre sequence \eqref{defsigma} induces a short exact sequence of $k_*$-modules:
\begin{equation}\label{defalgsigma}
\xymatrix{0 \lra k_* \ar[rr]_{k^R_*(\eta_F)} &&
k^R_*(F)\ar@/_1em/@{-->}[ll]_{\psi_*} \ar[rr]_-{k^R_*(\beta)}
&&k^R_*(\Sigma I) \lra 0}
\end{equation}
The retraction $\psi_*$, which is easily seen to be a
$k_*$-homomorphism, induces a $k_*$-linear section $\sigma_* \:
k^R_*(\Sigma I)\to k^R_*(F)$, which is natural in $F$ and $k$ as
well. So there is a natural isomorphism of $k_*$-modules
\begin{equation*}
k^R_*(F) \cong k_* \oplus k^R_*(\Sigma I),
\end{equation*}
given by $b \mapsto (\psi_*(b),k^R_*(\beta)(b))$, with inverse
$(c,a) \mapsto k^R_*(\eta_F)(c) + \sigma_*(a)$.

We define $\phi^k_F$ to be the composite
\[
\phi^k_F\: I/I^2[1] \cong F_*\otimes I[1] \xra{\pi_*\otimes 1} k_*\otimes I[1]
\xra{\zeta} (k\smash \Sigma I)_* \xra{\sigma_*} k_*^R(F),
\]
where $\zeta$ is the map as considered in \eqref{gradedcommutative}. Observe that
$\phi^k_F$ is a homomorphism of $F_*$-modules.

\begin{defn}\label{phigen}
We call $\phi^k_F$ the characteristic homomorphism of the admissible pair $(F, k, \pi)$.
If $k$ and $F$ are understood, we just write $\phi$.
\end{defn}
For another description of
$\phi$ based on a K\"unneth
spectral sequence compare
Remark \ref{ss}.

We defer the proof of the following fact to Section \ref{bilbasic}:

\begin{prop}\label{phiindep}
The characteristic homomorphism $\phi_F^k$ does not depend on the pro\-ducts on $F$ and
$k$.
\end{prop}

Recall that any regular quotient module $F$ of $R$ can be realized
as a regular quotient ring (Proposition \ref{realizereg}). The
following definition is meaningful by Proposition \ref{phiindep}:

\begin{defn}\label{charf}
The characteristic homomorphism of a regular quotient module $F=R/I$ is the
characteristic homomorphism $\phi^F_F\: I/I^2[1]\to F_*^R(F)$ of the multiplicative
admissible pair $(F, F, 1_F)$, where $F$ is endowed with an arbitrary product. We will
denote it by $\phi_F$ or simply by $\phi$ if $F$ is understood.
\end{defn}

\subsection{The characteristic bilinear form}\label{bilconstruction}

Assume that $(F, k, \pi)$ is an admissible pair. For brevity, we write ${}^{\pi}\phi$ for
the composition of the following $k_*$-homomorphisms
\[
{}^{\pi}\phi\: k_* \otimes_{F_*} I/I^2[1] \xra{1 \otimes \phi} k_* \otimes_{k_*} k_*^R(F)
\cong k_*^R(F).
\]
Recall the $R$-module map $\psi\: k\smash F\to k$ from Section \ref{phi} and the algebra
structure on $k_*^R(F)$ from Section \ref{notations}.

We define $b^k_F$ to be the composite of $k_*$-homomorphisms
\[
b^k_F\: (k_* \otimes_{F_*} I/I^2[1])^{\otimes 2}
\xra{{}^{\pi}\phi^{\otimes 2}} k^R_*(F)^{\otimes 2} \xra{m^k_F}
k_*^R(F) \xra{\psi_*} k_*.
\]

\begin{defn}
We call $b^k_F$ the characteristic bilinear form associated to the admissible pair
$(F,k,\pi)$.
\end{defn}

Observe that $b^k_F$ preserves the gradings and is natural in $F$
and $k$.

In the following, we write $\bar x$ for either of the elements
$(x+I^2)[1]\in I/I^2[1]$ or $1\otimes (x+I^2)[1]\in k_*
\otimes_{F_*} I/I^2[1]$ associated to some $x\in I$. The context
will make it clear which element is meant.

Associated to $b^k_F$ is the quadratic form $q^k_F\: k_*
\otimes_{F_*} I/I^2[1]\to k_*$, defined by $q^k_F(\bar x) =
b^k_F(\bar x\otimes \bar x)$ for $x\in I$. Note that $q^k_F$
doubles the degrees.

In the special situation of Example \ref{admcoupleFF} ($k=F$ as
$R$-modules and $\pi=1_F$), we write $b^{\nu}_\mu$ and
$q^{\nu}_\mu$ instead of $b_F^F$ and $q_F^F$, to keep track of the
products. If $\mu=\nu$, we simply write $b_F$ and $q_F$. If
$\mu=\nu^{\op}$, we write $b_{F^{\op}}^F$ and $q_{F^{\op}}^F$.
Whenever no confusion can be caused, we simply write $b$ and $q$.

If $(F,k,\pi)$ is multiplicative, its bilinear form $b^k_F$ is
determined by $b_F$ as well as by $b_k$. To describe the
relationship, let $k_*\otimes b_F$ denote the bilinear form on the
$k_*$-module $k_*\otimes_{F_*} I/I^2[1]$ determined by
\[
(k_*\otimes b_F)((1\otimes\bar x) \otimes (1\otimes\bar y)) = \pi_*(b_F(\bar x\otimes
\bar y)).
\]
Let moreover $\pi^*(b_k)$ be the bilinear form on $k_*\otimes_{F_*} I/I^2[1]$ determined
by
\[
\pi^*(b_k)((1\otimes\bar x) \otimes (1\otimes\bar y)) = b_k(\bar\pi_*(\bar x)\otimes
\bar\pi_*(\bar y)),
\]
where $\bar\pi_*\: I/I^2\to J/J^2$ is the canonical homomorphism and where $k=R/J$.

\begin{prop}\label{bilmultcouple}
The characteristic bilinear form of a multiplicative admissible pair $(F,k,\pi)$ is given
by $b^k_F = k_* \otimes b_F = \pi^*(b_k)$.
\end{prop}

\begin{proof}
This follows from naturality, by considering the admissible pairs $(F,F, 1_F)$, $(F, k,
\pi)$ and $(k,k, 1_k)$.
\end{proof}

The bilinear form $b^k_F$ will be determined for various $k$ and
$F$ in the next subsection. At this point, we can offer the
following general statement, which will be useful in the sequel.

\begin{prop}\label{oppositetrivial}
Let $(F,k,\pi)$ be a multiplicative admissible pair. Then the
characteristic bilinear form $b^k_{F^\op}$ of the admissible pair
$(F^\op, k, \pi)$ is trivial. In particular, $b^F_{F^\op}=0$ for a
quotient ring $F$.
\end{prop}

\begin{proof}
We first show that $b^F_{F^\op}=0$ for a quotient ring $F$. The natural left and right
actions of $F$ on $F\smash F$ and $F$ induce left actions of $F\smash F^\op$. The product
map $\mu\: F\smash F\to F$ respects these actions, and so $\mu_*\: F_*^R(F)\to F_*$ is a
map of left $F_*^R(F^\op)$-modules.
%Observe that the product on $F$ induces a left $F\smash F^\op$-module structure on both
%$F\smash F$ and $F$,\marginpar{I am still not very happy } moreover the product map
%$\mu\: F\smash F\to F$ is a map of left $F\smash F^\op$-modules. Hence $\mu_*\: F_*^R(F)
%\to F_*$ is left $F_*^R(F^\op)$-linear.
On $F_*^R(F^\op)$, the $F_*^R(F^\op)$-action is the same as the
one given by left multiplication in the algebra $F_*^R(F^\op)$. As
a consequence, we have for any $x,y\in I$ with residue classes
$\bar x, \bar y\in k_* \otimes_{F_*} I/I^2[1]$ (where $\cdot^\op$
denotes the product in $F_*^R(F^\op)$):
\[
\begin{tabular}{lll}
$b^F_{F^\op}(\bar x\otimes \bar y) $& $=\psi_*(\phi(\bar x)\cdot^\op \phi(\bar y)) =
\mu_*(\phi(\bar x) \cdot^\op \phi(\bar y))$
\\
&$ = \phi(\bar x)\cdot \mu_*(\phi(\bar y)) = 0,$
\end{tabular}
\]
because $\psi_*=\mu_*$ (second equality), $\mu_*$ is
$F_*^R(F^\op)$-linear (third equality) and $F_*$ is concentrated
in even degrees (fourth equality).

The statement for arbitrary multiplicative admissible pairs
$(F,k,\pi)$ now follows directly from Proposition
\ref{bilmultcouple}.
\end{proof}

\begin{cor}\label{bilcomm}
For a commutative quotient ring $F$, we have $b_F=0$.
\end{cor}

\subsection{The test case $F=R/x$}\label{bilbasic}

Assume that $(R/x, k, \pi)$ is an admissible pair, where $x\in
R_d$. We will first determine its characteristic homomorphism and
bilinear form.

We need some preparations. Applying $k^R_*(-)$ to the cofibre sequence  \eqref{cofx}
gives the short exact sequence of $k_*$-modules
\begin{equation}\label{ses_1gen}
0 \lra k^R_*(R) \xra{k^R_*(\eta_x)} k^R_*(R/x)
\xra{k^R_*(\beta_x)} k^R_*(\Sigma^{d+1}R) \lra 0.
\end{equation}
Because of $k_\odd=0$, $k^R_*(\Sigma^{d+1}R)\cong k_*[d+1]$ and
because $d$ is even, there exists a unique class $a_x \in
k^R_*(R/x)$ with $k^R_*(\beta_x)_*(a_x) =1[d+1]$. Therefore
\begin{equation}\label{smalliso}
k^R_*(R/x) \cong k_* \oplus k_*[d+1],
\end{equation}
where $1 \in k^R_*(R/x)$ corresponds to $(1,0)$ and $a_x \in k^R_*(R/x)$ to $(0,1[d+1])$.

\begin{rem}\label{freehomology}
By \eqref{smalliso}, the $k_*$-module $k^R_*(R/x)$ is $k_*$-free.
As a consequence,
$k^R_*(F)$ is $k_*$-free for any regular quotient $F$. Namely,
by Remark \ref{htyregquot}
and a K\"unneth isomorphism,
$k^R_*(F) \cong  \colim_k k^R_*(R/x_1)\otimes  \cdots
\otimes k^R_*(R/x_k)$ is $k_*$-free.
For another argument
based on
a K\"unneth
spectral sequence, see
Remark \ref{ss}.
\end{rem}

The $k_*$-module $k_* \otimes_{F_*} I/I^2[1]$ is freely generated
by $\bar x$. Therefore, $b=b^k_{R/x}$ and $q=q^k_{R/x}$ are
determined by the single element $b(\bar x \otimes \bar x) =
q(\bar x)$.

\begin{lem}\label{admcoll}
We have $\phi^k_{R/x}(\bar x) = a_x$ and $q^k_{R/x}(\bar x)\cdot 1 = a_x^2$.
\end{lem}

\begin{proof}
The first equality is a direct consequence of the definition of
$\phi$. For the second one, notice that by definition of $q$ and
by the first equality, we have
\[
q(\bar x)= \psi_*(k_*^R(\mu)(\kappa(a_x\otimes a_x))) = \psi_*(m^k_F(a_x\otimes
a_x))=\psi_*(a_x\cdot a_x).
\]
This implies the statement for dimensional reasons.
\end{proof}

We can now prove Proposition \ref{phiindep}:

\begin{proof}[Proof of Proposition \ref{phiindep}]
Observe first that $\phi^k_F$ is obviously independent on the product on $F$, since the
latter does not enter into its definition.

To show independence on $\mu_k$, we let $x\in I$ be arbitrary and show that
$\phi^k_F(\bar x)$ can be expressed without reference to $\mu_k$. Let $\bar \eta_F \: R/x
\to F $ be the unique factorization of $\eta_F\: R \to F$. Choose a product on $R/x$ such
that $\bar\eta_F$ is multiplicative (Proposition \ref{natargument}). Then the pair
$(\bar\eta_F, 1_F)$ is a morphism between the admissible pairs $(R/x, k, \pi\bar\eta_F)$
and $(F, k, \pi)$. Therefore, by naturality of the characteristic homomorphism, the
following diagram commutes:
\[
\xymatrix{(x)/(x^2)[1] \ar[d] \ar[r]^-{\phi^{k}_{R/x}} & k_*^R(R/x) \ar[d]
\\
I/I^2[1] \ar[r]^-{\phi^k_F} & k_*^R(F).}
\]
Now $\phi^{k}_{R/x}(\bar x)=a_x$ by Lemma \ref{admcoll}, which is defined independently
of the product on $k$. Hence so is $\phi^k_F(\bar x)$.
\end{proof}

We now aim to relate $q_{R/x}$ to Strickland's commutativity
obstruction $c(\mu_{R/x})$ (Proposition \ref{obstruction}).

\begin{prop}\label{identifyobstruction}
For a regular quotient ring $F=R/x$ with product $\mu$, we have $q_F(\bar x)=-c(\mu)\in
R_*/x$.
\end{prop}

\begin{proof}
The quadratic form $q_F$ on $(x)/(x^2)[1]\cong R_*/x\cdot \bar x$ on the one hand is
determined by $q=q_F(\bar x)=\mu_*(a_x^2)$ (we are in the situation where $\psi=\mu$).
The obstruction $c=c(\mu)$ on the other hand is characterized by the identity
$c(\beta\smash\beta) = \mu-\mu\tau$, where $\beta=\beta_x: F \to \Sigma^{|x|+1} R$ is
taken from the cofibre sequence \eqref{cofx}.
%and $\tau\: F\smash F\to F\smash F$ is the switch map.
Therefore we need to show that the maps $f_1=-q(\beta\smash\beta)$
and $f_2=\mu-\mu\tau$ coincide. We prove this using the
isomorphism of $F_*$-modules
\begin{equation}\label{geoalg}
d\kappa^*\: F^*_R (F\smash F) \lra \Hom_{R_*}(F^R_*(F)\otimes
F_*^R(F), F_*)
\end{equation}
given by composing the duality isomorphism $d$ from Section \ref{notations} with the one
induced by the K\"unneth isomorphism $\kappa=\kappa_\mu$.

First consider $(d\kappa^*)(f_1)$. Observe that by definition %(see \eqref{smalliso})
$k^R_*(\beta)(1)=0$ and $k^R_*(\beta)(a_x)=1$. From this, we easily deduce that
\[
(d\kappa^*)(f_1)(1 \otimes 1) = (d\kappa^*)(f_1)(a_x \otimes 1)=
(d\kappa^*)(f_1)(1 \otimes a_x) = 0
\]
and that
\[
(d\kappa^*)(f_1)(a_x \otimes a_x)=-q(k^R_*(\beta)\otimes
k^R_*(\beta))(a_x \otimes a_x)=q
\]
A sign is arising here according to \eqref{gradedcommutative}, because we let commute an
odd degree map, $k^R_*(\beta)$, with an odd degree element, $a_x$.

Now consider $(d\kappa^*)(f_2).$ As both $\mu$ and
$\tau\mu=\mu^\op$ are products on $F$, we have
\[
0=(d\kappa^*)(f_2)(1 \otimes 1) =(d\kappa^*)(f_2)(a_x \otimes
1)=(d\kappa^*)(f_2)(1 \otimes a_x).
\]
By definition of $q$, we have $(d\kappa^*)(\mu)(a_x \otimes
a_x)=\mu_*(a_x \cdot a_x)=q$ and moreover, as $a_x \cdot ^\op
a_x=0 \in F^R_*(F^\op)$ by Lemma \ref{admcoll} and Proposition
\ref{oppositetrivial},
\[
(d\kappa^*)(\mu^\op)(a_x \otimes a_x)=\mu_*(a_x \cdot^\op a_x)=0.
\]
It follows that $(d\kappa^*)(f_2)=(d\kappa^*)(f_1)$, which
concludes the proof.
\end{proof}

\subsection{Diagonal ring spectra}

The main aim of this subsection is to determine the characteristic
bilinear form of a diagonal regular quotient ring. More generally,
we consider $R$-rings $F$ which are obtained by smashing together
an arbitrary family of quotient $R$-ring spectra $F_i$. We specify
conditions on the $F_i$ which imply that $F$ is a quotient ring
and that the characteristic bilinear form $b_F$ is determined by
those of the $F_i$.

Suppose that $(F_i, \mu_i, \eta_i)_{i\geq 1}$ is a family of
$R$-ring spectra. There is an obvious way to endow a finite smash
product $F_1\smash\cdots\smash F_n$ with a product structure, by
mimicking the construction of the tensor product of finitely many
algebras. We now show that this construction extends to infinitely
many smash factors. Let $F=F_1\smash F_2\smash\cdots$ and let
$j_i\: F_i\to F$ be the natural maps.  The following statement
generalizes \cite{strickland}*{Prop.\@ 4.8}:

\begin{prop}\label{diagonalprop}
There is a unique $R$-ring structure on $F$ such that $j_k$
commutes with $j_l$ if $k\neq l$.
\end{prop}

\begin{proof}
There is an obvious right action of $F_n$ on $F(n)=
F_1\smash\cdots\smash F_n$. It extends in an evident way to
compatible $F_n$-actions on $F(i)$ for all $i\geq n$, which induce
an action $\psi_n\: F\smash F_n \to F$. We claim that the natural
maps $\pi_{n}\: [B\smash F_{n}, F]\to [B, F]$ induced by the units
$\eta_n\: R\to F_n$ are surjective for any $R$-module $B$. In
fact, we obtain a section of $\pi_{n}$ by associating to a map
$\alpha\: B\to F$ the composition $\psi_n(\alpha\smash 1)\:
B\smash F_n \to F$, because the diagram
\[
\xymatrix{B\smash F_n \ar[r]^-{\alpha\smash 1} & F\smash F_n
\ar[r]^-{\psi_n} & F
\\
B\ar[u]^-{1\smash\eta_n} \ar[r]^-\alpha & F
\ar[u]^-{1\smash\eta_n} \ar@{=}[ru]}
\]
commutes. As a consequence, we find that $[F^{\smash r}, F] \cong
\lim_n [F(n)^{\smash r}, F]$ for $r\geq 1$, by Milnor's exact
sequence. For the rest of the argument, we follow the proof of
\cite{strickland}*{Prop.\@ 4.8}.
\end{proof}

\begin{defn}\label{diagonal2}
We call $F$ with the product from Proposition \ref{diagonalprop}
the {\em smash ring spectrum}\/ of the $F_i$.
\end{defn}

Suppose now that $(F_i=R/I_i, \mu_i, \eta_i)_{i\geq 1}$ is a
family of quotient rings. Let $(F, \mu, \eta)$ be the smash ring
spectrum of the $F_i$ (Definition \ref{diagonal2}) and let
$I=I_1+I_2+\cdots$. We aim to express $b_F$ in terms of the
$b_{F_i}$ under conditions on the ideals $I_i$ which guarantee
that $F_*\cong R_*/I$ and that
\[
I/I^2 \cong \bigoplus_i R_*/I\otimes_{R_*} I_i/I_i^2.
\]
To begin with, note that the canonical homomorphisms
\[
R_*/(I_1+\cdots + I_k) \cong (F_1)_*\otimes\cdots\otimes (F_k)_*
\lra (F_1\smash\ldots \smash F_k)_*
\]
induce on passing to colimits a map $\theta\:
R_*/I=R_*/(I_1+I_2+\cdots) \to F_*$. Consider the following
hypotheses:
\begin{enumerate}
\item \label{theta} $\theta$ is an isomorphism; \item
\label{prodinter} $(I_1+ \cdots + I_{k-1}) \cdot I_k = (I_1 +
\cdots + I_{k-1})\cap I_k$ for all $k>1$.
\end{enumerate}

\begin{rem}
It may be interesting to note that in the case where $I_k=(x_k)$ for all $k$, hypothesis
(ii) is equivalent to the condition that $(x_1, x_2, \ldots)$ is a regular sequence. This
is easy to verify. The assumption that $R_*$ is a domain is essential here.
\end{rem}

\begin{prop}
Hypotheses (i) and (ii) are both satisfied if for $k>1$
\[
\Tor_{i,*}^{R_*}(R_*/(I_1+\cdots+I_{k-1}), R_*/I_k)=0\quad \forall
i>0.
\]
In particular, (i) and (ii) hold if $I_k$ is generated by a
sequence which is regular on $R_*/(I_1+\cdots+I_{k-1})$, for all
$k>1$.
\end{prop}

\begin{proof}
To show (i), we prove by induction that
\[
R_*/(I_1+\cdots + I_k) \cong (F_1\smash\cdots\smash F_k)_*.
\]
For the inductive step, it suffices to consider the K\"unneth spectral sequence
\[
E_{*,*}^2 = \Tor_{*,*}^{R_*}((F_1\smash\cdots \smash F_{k-1})_*,
(F_k)_*) \quad \Lra \quad (F_1\smash\cdots\smash F_k)_*,
\]
(see \cite{ekmm}*{IV.4}), which degenerates by assumption.

For (ii), recall that for ideals $J,K\subseteq R_*$, we have \cite{eisenbud}*{Exercise
A3.17}
\[
\Tor_{1,*}^{R_*}(R_*/J, R_*/K) = (J\cap K)/(J \cdot K).
\]

The last statement can be easily verified by using Koszul
complexes.
\end{proof}

The following fact must be well known. For lack of a reference, we
indicate its proof in Appendix \ref{commalg}.

\begin{prop}
Suppose that (ii) is satisfied. Then there is a canonical
isomorphism of $R_*/I$-modules
\begin{equation}\label{decompdiag}
I/I^2\cong \bigoplus_{i\geq 1}  R_*/I \otimes_{R_*} I_i/I_i^2.
\end{equation}
\end{prop}

We record the following immediate, well-known consequence:

\begin{cor}\label{imodi2}
Let $I \subset R_*$ be an ideal generated by a regular sequence
$(x_1, x_2, \ldots)$ and let $\bar x_i\in I/I^2$ denote the
residue classes of the $x_i$. Then there is an isomorphism of
$R_*/I$-modules $I/I^2\cong \bigoplus_{i} R_*/I\, \bar x_i$.
\end{cor}

The next proposition describes the characteristic bilinear form
associated to a smash ring spectrum. For the definition of the
bilinear forms $F_*\otimes b_{F_i}$ see the paragraph preceding
Proposition \ref{bilmultcouple}.

\begin{prop}\label{bildiagonal}
Let $F$ be the smash ring spectrum of quotient rings $F_i$ and
suppose that conditions (i) and (ii) above are satisfied. Then the
bilinear form $b_F$ is isomorphic to the direct sum of the
$F_*\otimes b_{F_i}$.
\end{prop}

\begin{proof}
Let $V_i=I_i/I_i^2[1]$, $V=I/I^2[1]$ and let $j_i\: F_i\to F$ be
the natural maps. As a consequence of naturality, the diagonal
terms of the bilinear form $b_F$ with respect to the decomposition
in condition (ii) are given by $F_*\otimes b_{F_i}$. Hence we need
to show that the off-diagonal terms of $b_F$ vanish. More
precisely, we must have $b_F(\bar x_k\otimes \bar x_l)=0$ for
$k\neq l$, $x_k\in I_k$ and $x_l\in I_l$. By definition, this
means that the composition
\begin{equation}\label{somecomp}
V_k \otimes V_l \to V\otimes V \xra{\phi_F\otimes \phi_F}
F_*^R(F)\otimes F_*^R(F) \xra{m^F_F} F_*^R(F) \xra{\mu_*} F_*
\end{equation}
has to be trivial, where the first map is induced by the
inclusions of $I_k$ and $I_l$ into $I$. By naturality, the
composition of the first two morphisms of \eqref{somecomp}
coincides with
\begin{equation*}
V_k  \otimes V_l \xra{\phi_{F_k}^F\otimes\phi_{F_l}^F}
F^R_*(F_k)\otimes F_*^R(F_l) \xra{F^R_*(j_k)\otimes F^R_*(j_l)}
F_*^R(F)\otimes F_*^R(F).
\end{equation*}
Because $j_k\: F_k\to F$ and $j_l\: F_l\to F$ commute, the composition of the last two
morphisms of \eqref{somecomp} with $F^R_*(j_k)\otimes F^R_*(j_l)$ coincides with
\begin{equation*}
F^R_*(F_k)\otimes F_*^R(F_l) \xra{F^R_*(j_k)\otimes F^R_*(j_l)} F_*^R(F)\otimes F_*^R(F)
\xra{m^F_{F^\op}} F_*^R(F) \xra{\mu_*} F_*.
\end{equation*}
Note that $m^F_{F^\op}$ can be viewed as the left action map of $F_*^R(F^\op)$ on itself
which is induced by the left action of $F\smash F^\op$ on itself. Now $\mu_*\: F_*^R(F)
\to F_*$ is left $F_*^R(F^\op)$-linear, as we have noted earlier. Because $V_k$ and $V_l$
are concentrated in odd degrees, an argument as in the proof of Proposition
\ref{oppositetrivial} shows that \eqref{somecomp} is zero.
\end{proof}

We close this section by determining the characteristic bilinear
form $b_F$ of a diagonal regular quotient ring $F$.

\begin{prop}\label{bildiag}
Let $(x_1, x_2, \ldots)$ be a regular sequence in $R_*$ generating an ideal $I\subset
R_*$. Suppose that $\mu_i$ are products on $R/x_i$ and let $F=R/I=R/x_1\smash
R/x_2\smash\cdots$ be the induced diagonal regular quotient ring. Then the characteristic
bilinear form $b_F\: I/I^2[1]\otimes_{F_*} I/I^2[1]\to F_*$ is diagonal with respect to
the basis $\bar x_1, \bar x_2, \ldots$ and $b_F(\bar x_i\otimes \bar x_i)\equiv -c(\mu_i)
\mod I$.
\end{prop}

\begin{proof}
Combine Propositions \ref{identifyobstruction} and
\ref{bildiagonal}.
\end{proof}

\section{The homology algebra}\label{homology}

The aim of this section is to study the homology algebra $k_*^R(F)$ for an admissible
pair $(F, k, \pi)$, with its natural product $m^k_F$ from Section \ref{bilconstruction}.

\subsection{The main result and some consequences}\label{algstructure}

Before stating the main result, we recall the definition and the universal property of
Clifford algebras.

Let $M_*$ be a graded quadratic module, \ie\ a graded module over
a graded commutative ring $k_*$, endowed with a quadratic form
$q\: M_* \to k_*$ which doubles degrees (for instance the
quadratic form associated to a degree-preserving bilinear form).
Let $\T(M_*)$ denote the tensor algebra over $k_*$, with its
natural grading. The Clifford algebra $\Cl(M_*, q)$ is defined as
\[
\Cl(M_*, q) = \T(M_*)/(x\otimes x-q(x)\cdot 1; \ x\in M_*).
\]
As the ideal $(x\otimes x-q(x)\cdot 1; \ x\in M_*)$ is homogenous,
$\Cl(M_*, q)$ inherits a grading from $\T(M_*)$. Up to unique
isomorphism, $\Cl(M_*, q)$ is characterized by the following
universal property: Any degree-preserving $k_*$-linear map $f\:
M_* \to A_*$ into a graded $k_*$-algebra $A_*$ such that
$f(x)^2=q(x)\cdot 1$ for all $x\in M_*$ lifts to a unique algebra
map $\Cl(M_*, q) \to A_*$.

\begin{thm}\label{algebra_main}
Let $(F=R/I,k, \pi)$ be an admissible pair. Then the
characteristic homomorphism
\[
{}^\pi\phi\: k_*\otimes_{F_*} I/I^2[1] \lra k^R_*(F)
\]
lifts to a natural homomorphism of $k_*$-algebras
\[
\Phi\:  \Cl(k_*\otimes_{F_*} I/I^2[1], q^k_F) \lra k^R_*(F).
\]
If $F$ is a regular quotient, then $\Phi$ is an isomorphism.
\end{thm}

We will prove this result in Section \ref{mainproof} and draw some
consequences now. Let us first spell out the following important
special cases:

\begin{cor}\label{coralgiso}
Let $F=R/I$ be a regular quotient ring. Then there is a natural $F_*$-algebra isomorphism
\[
F_*^R(F)\cong \Cl(I/I^2[1], q_F).
\]
\end{cor}

\begin{cor}\label{algffop}
Let $F=R/I$ be a regular quotient ring. Then there is an
$F_*$-algebra isomorphism
\[
F_*^R(F^\op)\cong\Lambda(I/I^2[1]).
\]
Under this isomorphism, the homomorphism $(\mu_F)_*\:
F_*^R(F^\op)\to F_*$ corresponds to the canonical augmentation
$\epsilon\: \Lambda(I/I^2[1]) \to F_*$.
\end{cor}

\begin{proof}
The first statement follows from the fact that $q^F_{F^\op}=0$, by Proposition
\ref{oppositetrivial}. For the second statement, note that the map $(\mu_F)_*$ is
determined as the unique $F_*^R(F^\op)$-bilinear map which is trivial on the image of
$\phi$. The augmentation $\epsilon$, in turn, is a map of algebras, hence
$\Lambda(I/I^2[1])$-bilinear, and it is trivial on $I/I^2[1]$. Hence the two maps
coincide.
\end{proof}

\begin{rem}
Let $(F,k,\pi)$ be a multiplicative admissible pair, with $F=R/I$
a regular quotient ring. From Corollary \ref{algffop} and
Proposition \ref{bilmultcouple}, we deduce that there is an
isomorphism of $k_*$-algebras
\[
k_*^R(F^\op) \cong \Lambda(k_*\otimes_{F_*} I/I^2[1]).
\]
\end{rem}

We can be more explicit in the case of a regular quotient ring $F=R/I$ if we fix a
regular sequence $(x_1, x_2, \ldots)$ generating $I$. By Corollary \ref{imodi2}, this
choice determines an isomorphism $I/I^2\cong\bigoplus_{i} F_* \bar x_i$, where $\bar x_i$
denote the residue classes, as usual. Letting $a_i=\phi(\bar x_i)\in F_*^R(F)$, we have
\begin{equation}\label{ffopexplicit}
F_*^R(F^\op) \cong \Lambda(a_1, a_2, \ldots).
\end{equation}
Assume now that $F$ is diagonal and let $c_i\in R_*/x_i$ be the commutativity obstruction
of $R/x_i$ of Proposition \ref{obstruction} and let $\bar c_i$ be its residue class in
$F_*$. Using the explicit description of $b_F$ (and hence $q_F$) from Proposition
\ref{bildiag}, we find:
\begin{equation}\label{ffexplicit} F_*^R(F) \cong
\T(a_1, a_2, \ldots)/(a_i^2+\bar c_i \cdot 1,\, a_k a_l + a_l a_k; \,i\geq 1,\, k \neq
l).
\end{equation}

We add an example to illustrate the usefulness of the naturality of the isomorphism
$\Phi$ in Theorem \ref{algebra_main}.

\begin{exmp}\label{exa}
Let $R=\H\Z$ and $p$ be a prime. Recall that $R$, $F=\H\Z/p^4$ and $G=\H\Z/p^3$ are
commutative $\SS$-algebras and that the canonical map $F \to G$ corresponding to the
inclusion $I=(p^4) \to J=(p^3)$ is multiplicative \cite{ekmm}*{IV.2}. Multiplication by
$p^4$ and $p^3$ induces isomorphisms $\Z/p^4\cong I/I^2$ and $\Z/p^3\cong J/J^2$,
respectively. Under these identifications, the map $I/I^2\to J/J^2$ corresponds to $p\:
\Z/p^4\to\Z/p^3$. For any $(G, k, \pi)$ admissible, the map of $k_*$-algebras
$k_*^R(F)\to k_*^R(G)$ identifies with $\Lambda_{k_*}(a) \to \Lambda_{k_*}(b),$ $a
\mapsto p\cdot  b$. If $k=\H\Z/p$, this map is trivial, if $k=\H\Z/p^2$, it is
non-trivial.
\end{exmp}

\begin{rem}\label{ss}
\!\!\!\!\footnote{This remark has been suggested by the referee.}
The K\"unneth spectral sequence
(see \cite{ekmm}*{IV.4})
\begin{equation}\label{kss}
E^2_{p,q} = \Tor_{p,q}^{R_*}(k_*,F_*) \ \Longrightarrow
k_{p+q}^R(F)
\end{equation}
is a multiplicative spectral sequence of
$k_*$-algebras, see \cite{b-l}*{Lemma 1.3}. By standard
techniques,
\[
\Tor_{*,*}^{R_*}(k_*,F_*) \cong
\Lambda(k_*\otimes_{F_*} I/I^2[1])
\]
as $k_*$-algebras (this follows for
instance from \cite{ml}*{VII.6, Exercise 3}).
For dimensional reasons, the
elements of $k_*\otimes_{F_*} I/I^2[1]$
are permanent cycles and thus by
multiplicativity, the spectral
sequence collapses. As
$\Lambda(k_*\otimes_{F_*} I/I^2[1])$
is a free $k_*$-module, there are
no {\it additive} extensions and hence
$k_{*}^R(F)\cong \Lambda(k_*\otimes_{F_*} I/I^2[1])$
as $k_*$-modules.
The proof of Theorem \ref{algebra_main}
can be seen as resolving the
{\it multiplicative} extensions in the
spectral sequence.

The characteristic homomorphism
${}^{\pi}\phi\: k_* \otimes_{F_*} I/I^2[1]
\to k_*^R(F)$
can also be considered from the
point of view of this spectral
sequence. Let
$\mathscr{F}_{0}\subset\mathscr{F}_{1}
\subset \cdots
\subset k_*^R(F)$
be the filtration naturally
associated to the spectral
sequence.
Consider the short exact
sequence
\begin{equation}\label{defalgsigmar}
0 \lra \mathscr{F}_{0}\lra
\mathscr{F}_{1}\lra
E^\infty_{1,*}\lra 0.
\end{equation}
The retraction $\psi_*$
from \eqref{defalgsigma}
induces a natural retraction
in \eqref{defalgsigmar}.
Therefore we obtain
a natural isomorphism
$\mathscr{F}_{1}\cong \mathscr{F}_{0}
\oplus E^\infty_{1,*}$.
We can
show that the composition
\[
k_* \otimes_{F_*} I/I^2[1] \cong
\Tor_{1,*}^{R_*}(k_*,F_*)=
E^\infty_{1,*}\subset
\mathscr{F}_{1}\subset k_*^R(F)
\]
coincides with ${}^{\pi}\phi$.
\end{rem}

\subsection{Proof of Theorem \ref{algebra_main}}\label{mainproof}

To begin with, suppose that $F=R/x$, for some $x\in R_*$, and let $a_x=\phi^k_{R/x}(\bar
x)$. Then $k_*^R(F)\cong k_*\oplus k_* a_x$ by \eqref{smalliso} and $a_x^2 = q^k_x(\bar
x)\cdot 1$ by Lemma \ref{admcoll}. Hence $\Phi$ is an isomorphism of algebras
\begin{equation}\label{algrmodx}
\T_{k_*}(a_x)/(a_x^2-q^k_x(\bar x)\cdot 1) \cong k_*^R(R/x),
\end{equation}
which is exactly the statement of the theorem for $F=R/x$.

Assume now that $F$ is a quotient ring of $R$. By the universal property of Clifford
algebras, the lift $\Phi$ exists if and only if $\phi^k_F(\bar x)^2 =q^k_F(\bar x) \cdot
1$ for all $x \in I$. Fix $x\in I$ and consider the natural map $j\: R/x \to F$. There
exists a product $\mu_x$ on $R/x$ such that $j$ is multiplicative by Proposition
\ref{natargument}. Now the inclusion $(x) \subset I$ induces a commutative diagram of the
form
\[
\xymatrix{ (x)/(x)^2[1] \ar[r]^-{\phi^k_{R/x}} \ar[d] & k_*^R(R/x)
\ar[d]^-{k_*^R(j)}
\\
I/I^2[1]\ar[r]^-{\phi^k_{F}}  & k_*^R(F).}
\]
As $j$ is multiplicative, $k_*^R(j)$ is a map of algebras. We thus obtain
\[
\phi^k_F(\bar x)^2 = k_*^R(j)(\phi^k_{R/x}(\bar x)^2) =
k_*^R(j)(q^k_{R/x}(\bar x)\cdot 1) = q^k_F(\bar x)\cdot 1,
\]
by Lemma \ref{admcoll} and by naturality of $q$. It follows that $\phi^k_F$ lifts to an
algebra map $\Phi=\Phi^k_F$, as asserted.

Suppose now that $F$ is a regular quotient ring. To show that $\Phi$ is an isomorphism,
it suffices to prove this for the case where $I$ is generated by a finite regular
sequence $(x_1, \ldots, x_n)$. The general case then follows easily by passing to
colimits. Let $i_i\: (x_i)\to I$ denote the inclusions and $\bar i_i\: (x_i)/(x_i)^2\to
I/I^2$ the induced maps. As before, we choose products $\mu_i$ on the $R/x_i$ such that
the natural maps $j_i\: R/x_i\to F$ are multiplicative.
%Let $q^k_i$ be the characteristic quadratic forms associated to the $\mu_i$.
Consider the diagram of $k_*$-modules
\[
\xymatrix@C=5em{  \bigotimes\limits_{i=1}^n k_*\otimes \Cl((x_i)/(x_i)^2, q^k_{R/x_i})
\ar[r]^-{k_*\otimes\bigl(\bigotimes \Phi^k_{R/x_i}\bigr)} \ar[d]_-{\bigotimes (1\otimes
\bar i_i)} &
\bigotimes\limits_{i=1}^n k_*^R(R/x_i) \ar[d]^-{\bigotimes k_*(j_i)}\\
\bigotimes\limits_{i=1}^n k_*\otimes \Cl(I/I^2, q^k_F)
\ar[r]^-{k_*\otimes\bigl(\bigotimes \Phi^k_{F}\bigr)} \ar[d] &
\bigotimes\limits_{i=1}^n k_*^R(F) \ar[d]\\
k_*\otimes \Cl(I/I^2, q_F^k)\ar[r]^-{k_*\otimes \Phi^k_{F}}  & k_*^R(F).}
\]
The two lower vertical maps are given by the multiplication maps
of the respective algebras. The top square commutes because
$\Phi^k_F$ is natural in $F$ and the bottom one because
$\Phi^k_{F}$ is a morphism of algebras. The top horizontal map is
an isomorphism by \eqref{algrmodx}. As $I/I^2\cong
\bigoplus_{i=1}^n F_* \bar x_i$ by Corollary \ref{imodi2}, the
composite of the two left vertical maps is an isomorphism by
\cite{bourbaki}*{Ch.\@ VI, §9.3, Cor.\@ 3}. We easily check that
the composite of the two right vertical maps is just the K\"unneth
morphism and therefore an isomorphism. It follows that
$\Phi^k_{F}$ is an isomorphism, as asserted. \hfill \qed

\subsection{The antipode}\label{antipode}
Theorem \ref{algebra_main} allows us to give a neat description of the antipode (or
conjugation) homomorphism $\tau_*\: F_*^R(F)\to F_*^R(F)$ induced by the switch map
$\tau\: F\smash F\to F\smash F$. For this, we recall a definition from the theory of
Clifford algebras. Let $\Cl(M_*, q)$ be the Clifford algebra on a quadratic graded module
$M_*$. Then the {\em principal automorphism} $\alpha$ is the uniquely determined algebra
automorphism  of $\Cl(M_*, q)$ whose restriction to $M_*$ is given by $\alpha(m)=
(-1)^{|m|}m $.

\begin{prop}\label{switch}
Let $F$ be a regular quotient ring. Under the isomorphism from Corollary \ref{coralgiso},
the morphism $\tau_*\: F_*^R(F) \to F_*^R(F)$ corresponds to the principal automorphism
\[
\alpha\: \Cl(I/I^2[1], q_F) \to \Cl(I/I^2[1], q_F).
\]
\end{prop}

\begin{proof}
Since the switch map
$\tau\: F\smash F\to F\smash F$ is a ring isomorphism,
$\tau_* \: F_*^R(F) \to F_*^R(F)$
is an algebra isomorphism. It therefore suffices to check that
$\tau_*(\phi(\bar x))= - \phi(\bar x)$ for $x\in I$. Because there
is always  a product on $R/x$ such that the natural map $R/x\to F$
is multiplicative (Proposition \ref{natargument}) and because
$\phi$ is natural, we may therefore restrict to the case where
$F=R/x$. We set $d=|x|$. Recall that $d$ is even.

Let $a_x=\phi(\bar x)$. Then $F_*^R(F) = F_* 1\oplus F_* a_x$ by \eqref{smalliso}.
Clearly, we have $\tau_*(1)=1$. We therefore need to show that $\tau_*(a_x)= -a_x$.

We prove this by considering the canonical homomorphism
\[
\iota\: F_*^R(F) \cong (R\smash F\smash F)_* \xra{(\eta\smash 1 \smash 1)_*} (F\smash
F\smash F)_* = F_*^R(F\smash F)
\]
from the homotopy groups of $F\smash F$ to its homology groups. As $\iota$ is injective
($\mu$ induces a retraction), it suffices to prove that $F^R_*(\tau)(\iota(a_x)) =
-\iota(a_x)$. We do this by first identifying $\iota(a_x)$ and then computing
$F^R_*(\tau)(\iota(a_x))$.

To simplify the notation, we identify $F_*^R(F\smash F)$ with $F_*^R(F) \otimes_{F_*}
F_*^R(F)$ via the K\"unneth isomorphism and $(R\smash M)_*$, as well as $(M\smash R)_*$,
with $M_*$, for any $R$-module $M$.

To determine $\iota(a_x)$, we start by noting that for dimensional
reasons and as
$F_{\odd}=0$, we have
\[
\iota _*(a_x)= r \cdot  1 \otimes a_x + s \cdot a_x \otimes 1,
\]
where $r,s \in F_0$. Consider the commutative diagram
\begin{equation*}
\begin{array}{c}
\xymatrix{ R\smash F\smash F  \ar[d]_-{1 \smash 1\smash \beta} \ar[r]^-{\eta\smash
1\smash 1} & F\smash F\smash F  \ar[d]^-{1 \smash 1\smash \beta}
\\
R\smash F\smash R  \ar[r]^-{\eta \smash 1\smash 1} & F\smash F\smash R .}
\end{array}
\end{equation*}
The composition of the upper and the right morphisms induces
\begin{align*}
F_*^R(1\smash \beta)(\iota _*(a_x)) & = 1 \otimes F_*^R(\beta)(r \cdot  1 \otimes a_x + s
\cdot a_x \otimes 1)
\\
& = r \cdot F_*^R( \beta)(a_x)= r \cdot 1,
\end{align*}
whereas the composition along the two other edges of the diagram induces
$\eta_*(F_*^R(\beta)(a_x))= \eta_*(1)=1$. It follows that $r=1$.

The computation of $s$ requires another commutative diagram, namely
\begin{equation*}
\begin{array}{c}
\xymatrix{ R\smash F\smash F  \ar[d]_-{1 \smash \mu} \ar[r]^-{\eta\smash 1\smash 1} &
F\smash F\smash F \ar[d]^-{1 \smash \mu}
\\
R\smash F \ar[r]^-{\eta \smash 1} & F\smash F .}
\end{array}
\end{equation*}
With the same strategy as above, we obtain that
\[
F_*^R(\mu)(\iota _*(a_x))= F_*^R(\mu)(1 \otimes a_x + s \cdot a_x \otimes 1)= a_x + s
\cdot a_x
\]
and $\eta_*(\mu_*(a_x))=0$. This implies that $s=-1$.

We now consider the two maps $i=1\smash\eta$, $j=\eta\smash 1 \: F \to F\smash F$. The
induced morphisms in homology satisfy $F_*^R(i)(a_x)=a_x \otimes 1$ and $F_*^R(j)(a_x)=1
\otimes a_x$, respectively. Since $\tau \circ i = j$ and $\tau^2=1_{F\smash F}$, we
deduce that
\begin{equation}\label{tauaction}
F_*^R(\tau)(a_x \otimes 1)=1 \otimes a_x, \; F_*^R(\tau)(1 \otimes a_x)=a_x \otimes 1.
\end{equation}
Therefore, we have shown that
\[
F^R_*(\tau)(\iota(a_x) = F_*(\tau)(1\otimes a_x - a_x \otimes 1) = a_x\otimes 1 - 1
\otimes a_x = - \iota(a_x),
\]
which concludes the proof.
\end{proof}

\section{The cohomology algebra}\label{cohomology}

The aims of this section are to to give a natural description of
the cohomology module $k^*_R(F)$ for an admissible pair, to
identify the derivations $\theta\: F\to k$ in case the pair is
multiplicative and to identify canonically the cohomology algebra
$F^*_R(F)$ for a regular quotient ring $F$.

\subsection{The cohomology of admissible pairs}\label{basiccoho}

Let $(F, k, \pi)$ be an admissible pair. Using our identification of homology $k_*^R(F)$
and Kronecker duality, we derive an analogous expression for cohomology $k^*_R(F)$. As
for homology, we aim for an isomorphism which is natural in both $F$ and $k$.

For this, we need to modify the category of admissible pairs. We keep the objects, but
declare a morphism $(F, k, \pi) \to (G, l, \pi')$ in the new category to be a pair of
ring maps $f\: G\to F$ and $g\: k\to l$ such that
\begin{equation*}
\begin{array}{c}
\xymatrix{  F \ar[r]^-{\pi}  & k \ar[d]^-{g}
\\
G \ar[u]^-f \ar[r]^-{\pi'} & l}
\end{array}
\end{equation*}
commutes. We refer to this category as the {\em bivariant category
of admissible pairs} (this category
is also named the twisted
arrow category of the
category of admissible pairs). Similarly, we refer to the full subcategory
spanned by the multiplicative admissible pairs as the {\em
bivariant category of multiplicative admissible pairs}. We use the
expression ``natural in $F$ and $k$'' in an analogous sense as for
ordinary admissible pairs.

Cohomology $k^*_R(F)$ defines a functor on the bivariant category of admissible pairs.
Ordinary Kronecker duality $d\: k^*_R(F) \to D_{k_*}(k_*^R(F))$ (where $D_{k_*}(-) =
\Hom_{k_*}(-,k_*)$) is not appropriate to study this functor, because $D_{k_*}(k_*^R(F))$
is not functorial on the bivariant category. We can get around this inconvenience by
defining a modified version of Kronecker duality, which is of the form
\begin{equation}\label{kronmod}
d'\: k^*_R(F)\lra \Hom_{F_*}^*(F_*^R(F), k_*).
\end{equation}
It associates to a map $f\: F\to k$ the homomorphism
\[
F_*^R(F) \xra{F_*^R(f)} F_*^R(k) \xra{(\pi\smash 1)_*} k_*^R(k) \xra{(\mu_k)_*} k_*.
\]
We leave the easy verification of the fact that $d'$ is a natural
transformation between functors defined on the bivariant category
to the reader.

We will need to work with the profinite topology on $k^*_R(M)$ for $R$-modules $M$. This
is discussed in detail in \cite{sw}*{\S 2}, following ideas of \cite{boardman}. Recall
that for any graded $k_*$-module $N_*$, $D_{k_*}(N_*)$ carries a natural linear topology,
the dual-finite topology \cite{boardman}*{Def.\@ 4.8}, which is complete and Hausdorff.

We endow $\Hom_{F_*}(M_*, k_*)$, for a graded $F_*$-module $M_*$,
with the linear topology inherited
from the dual-finite topology on
$D_{k_*}(k_*\otimes_{F_*} M_*)$ under the adjunction isomorphism
\begin{equation}\label{adjunction}
\Hom^*_{F_*}(M_*, k_*) \cong D_{k_*}(k_*\otimes_{F_*} M_*).
\end{equation}
By a slight abuse of terminology, we refer to this topology as the
dual-finite topology, too. By naturality (in the variable $M_*$)
of \eqref{adjunction}, the function $M_*\mapsto\Hom^*_{F_*}(M_*,
k_*)$ gives rise to a functor from the category of $F_*$-modules
to the category of complete Hausdorff $k_*$-modules. As $d'$
agrees with the following composition (the unlabelled maps are the
canonical ones)
\begin{equation}\label{sthing}
k^*_R(F) \xra{d} D_{k_*}(k^R_*(F)) \lra D_{k_*}(k_*\otimes_{F_*} F_*^R(F)) \cong
\Hom^*_{F_*}(F_*^R(F), k_*),
\end{equation}
it is continuous, since the Kronecker homomorphism $d$ is
continuous.

\begin{prop}\label{cohomologyiso}
Let $(F=R/I,k, \pi)$ be an admissible pair. Then there exists a
natural continuous homomorphism of $k_*$-modules
\[
\Psi\: k^*_R(F) \lra \Hom_{F_*}^*(\Lambda(I/I^2[1]), k_*).
\]
If $F$ is a regular quotient ring, $\Psi$ is a homeomorphism.
\end{prop}

\begin{proof}
We define $\Psi$ as the composition of continuous homomorphisms
\begin{align*}
k^*_R(F) & \xra{d'} \Hom_{F_*}^*(F_*^R(F^\op), k_*) \xra{\Phi^*}
\Hom^*_{F_*}(\Lambda(I/I^2[1]), k_*),
\end{align*}
where $\Phi^*$ is induced by the homomorphism $\Phi\: \Lambda(I/I^2)\to F_*^R(F^\op)$
(Theorem \ref{algebra_main}, Proposition \ref{oppositetrivial}). If $F$ is a regular
quotient ring, the Kronecker duality homomorphism $d$ is a homeomorphism, as $k_*^R(F)$
is free over $k_*$ see \cite{sw}*{Prop.\@ 2.5}. Furthermore, it follows from Theorem
\ref{algebra_main} that the canonical homomorphism $k_*\otimes_{F_*} F_*^R(F) \to
k_*^R(F)$ is an isomorphism. Thus all the maps in \eqref{sthing} are homeomorphisms and
hence $d'$ as well. Moreover, $\Phi$ is an isomorphism by Theorem \ref{algebra_main} and
hence $\Phi^*$ is a homeomorphism by functoriality. It follows that $\Psi$ is an
homeomorphism.
\end{proof}

\subsection{Derivations of regular quotient rings}\label{detect}

Let $F$ be an $R$-ring and $M$ an $F$-bimodule. Recall that a map $\theta\: F\to \Sigma^i
M$ in $\DR$ is called a (homotopy) derivation if the diagram
\begin{equation}\label{htyderivation}
\begin{array}{c}
\xymatrix{ F\smash F \ar[rr]^-{1\smash \theta \vee \theta\smash 1} \ar[d]^-{\mu_F} &&
(F\smash \Sigma^i M) \vee (\Sigma^i M\smash F) \ar[d]^-{}
\\ F \ar[rr]^-{\theta} && \Sigma^i M}
\end{array}
\end{equation}
commutes, where the unlabelled map is induced by the left and
right actions of $F$ on $M$. We write $\DDer^i_R(F, M)$ for the
set of all such derivations and $\DDer^i_R(F)$ for $\DDer^i_R(F,
F)$.

Suppose that $(F=R/I,k, \pi)$ is a multiplicative admissible pair.
Then $k$ is an $F$-bimodule in a natural way, and so we may
consider $\DDer_R^*(F, k)$. We endow $\DDer_R^*(F,k)$ with the
subspace topology induced by the profinite topology on $k^*_R(F)$.

We now define a natural transformation
\[
\psi\: \DDer_R^*(F, k) \lra \Hom^*_{F_*}(I/I^2[1], k_*)
\]
between functors on the bivariant category of multiplicative
admissible pairs with values in the category of  topological
$k_*$-modules. We set $\psi$ to be the composition
\begin{equation*}
\DDer_R^*(F, k)  \subseteq k^*_R(F) \xra{\Psi} \Hom^*_{F_*}(\Lambda(I/I^2[1]), k_*)
 \xra{\iota^*} \Hom^*_{F_*}(I/I^2[1], k_*),
\end{equation*}
where $\Psi$ is the homomorphism from Proposition \ref{cohomologyiso} and where $\iota$
denotes the canonical injection $I/I^2[1] \to \Lambda(I/I^2[1])$.

\begin{prop}\label{derivations}
Suppose that $(F,k, \pi)$ is a multiplicative admissible pair and
that both $F=R/I$ and $k$ are regular quotient rings.
\begin{enumerate}\itemsep2pt
\item The homomorphism
\[\psi\: \DDer_R^*(F, k) \to \Hom^*_{F_*}(I/I^2[1], k_*)
\]
is a natural homeomorphism.
\item The composition
\[
\Hom^*_{F_*}(I/I^2[1], k_*) \xra{\psi^{-1}} \DDer_R^*(F,k) \subseteq k^*_R(F).
\]
is independent of the products on $F$ and $k$.
\end{enumerate}
\end{prop}

The proof of Proposition \ref{derivations} requires some
preparations and will be given at the end of this subsection.

To be able to detect derivations, we now relate homotopy
derivations with algebraic ones. We denote by $\Der_{k_*}^*(A_*,
M_*)$ the derivations from a $k_*$-algebra $A_*$ to an
$A_*$-bimodule $M_*$ and write $\Der^*_{k_*}(A_*)$ for
$\Der_{k_*}^*(A_*, A_*)$. The grading convention is that
$\partial\in\Der_{k_*}^*(A_*, M_*)$ satisfies
\[
\partial(a\cdot b) = \partial(a)\cdot b + (-1)^{|\partial|\cdot|a|} a\cdot \partial(b).
\]

\begin{lem}\label{derivationlemma}
Let $(F, k, \pi)$ be a multiplicative admissible pair, where
$F=R/I$ is a regular quotient, and let $\bar k$ be $k$, endowed
with a second (not necessarily different) product. Then the
Hurewicz homomorphism
\begin{equation}\label{hurewicz}
h = k_*^R(-)\: k^*_R(F) \lra \Hom^*_{k_*}(\bar k_*^R(F), \bar
k_*^R(k))
\end{equation}
restricts to a monomorphism
\[
\bar h\: \DDer_R^*(F, k) \lra \Der^*_{k_*}(\bar k_*^R(F), \bar k_*^R(k)).
\]
The induced commutative diagram
\[
\xymatrix{\DDer_R^*(F, k) \ar[r]^-{\bar h} \ar[d]_-{\incl} & \Der^*_{k_*}(\bar k_*^R(F),
\bar k_*^R(k)) \ar[d]^-{\incl}
\\
k^*_R(F) \ar[r]^-h & \Hom^*_{k_*}(\bar k_*^R(F), \bar k_*^R(k))}
\]
is a pullback diagram. Explicitly, this means that the derivations
are precisely those maps in $k^*_R(F)$ which induce derivations on
applying $\bar k_*^R(-)$.
\end{lem}

\begin{proof}
Applying the functor $h=\bar k_*^R(-)$ to the diagram \eqref{htyderivation} and
precomposing with the K\"unneth map $\kappa_{\bar k}\: \bar k_*^R(F)\otimes_{k_*} \bar
k_*^R(F) \to \bar k_*^R(F\smash F)$ shows that a derivation $\theta\: F\to \Sigma^i k$
induces a derivation $h(\theta)$ on the homology algebra $\bar k_*^R(F)$. Hence $h$,
which is monomorphic (see Section \ref{notations}), restricts to a monomorphism $\bar h$,
as asserted.

For the second statement, we need to verify, for $\theta\in
k^*_R(F)$, the equivalence
\begin{equation}\label{somequiv}
\theta\in\DDer_R^*(F, k)  \ \Longleftrightarrow \ h(\theta)\in
\Der^*_{k_*}(\bar k_*^R(F), \bar k_*^R(k)).
\end{equation}
We have shown ``$\Rightarrow$'' above and now prove
``$\Leftarrow$''. By definition of $\bar k_*^R(F)$ and $\bar
k_*^R(k)$, $h(\theta)$ is a derivation if the diagram obtained by
applying $\bar k_*^R(-)$ to \eqref{htyderivation} and precomposing
with $\kappa_{\nu_k}$ commutes. This implies that
\eqref{htyderivation} commutes (see Section \ref{notations}), \ie\
that $\theta$ is a derivation.
\end{proof}

Using Lemma \ref{derivationlemma}, we now construct certain
derivations in $\DDer_R^*(F)$. We first consider the case $F=R/x$.
Recall the maps $\beta_x, \eta_x$ from \eqref{cofx}. We refer to
the composition
\begin{equation}\label{bocksteinx}
Q_x\: R/x \xra{\beta_x} \Sigma^{d+1} R \xra{\eta_x} \Sigma^{d+1}
R/x
\end{equation}
as the {\em Bockstein operation}\/ associated to $x$.

The following lemma is already known from Strickland
\cite{strickland}. Let $ y \in
D_{R_*/(x)}\bigl((x)/(x^2)[1]\bigr)$ denote the dual of $\bar x
\in (x)/(x)^2[1]$.

\begin{lem}\label{bocksteinisderi}
The Bockstein operation $Q_x\: R/x \to \Sigma^{|x|+1} R/x$ is a
derivation for any product on $R/x$. It satisfies $\psi(Q_x)= y$.
\end{lem}

\begin{proof}
We have $(R/x^\op)_*^R(R/x) \cong \Lambda(a)$ with $a=\phi(\bar x)$, by Corollary
\ref{algffop}. Applying $(R/x^\op)_*^R(-)$ to the cofibre sequence \eqref{cofx}, we find
that under this isomorphism, $(R/x^\op)_*^R(Q_x)$ corresponds to
$\tfrac{\partial}{\partial a}\: \Lambda(a)\to \Lambda(a)$. Therefore, by Lemma
\ref{derivationlemma}, $Q_x$ is a derivation, with $\psi(Q_x)=y$.
\end{proof}

\begin{rem}\label{Qinduced}
The proof shows that $F^R_*(Q_x)$ corresponds to $\tfrac{\partial}{\partial a}$ under the
isomorphism $F_*^R(R/x)\cong \Lambda(a)$, where $a=\phi(\bar x)$.
\end{rem}

Next, we construct derivations in $\DDer_R^*(F)$ for an arbitrary regular quotient ring
$F=R/I$. Let $(x_1, x_2 ,\ldots)$ be a regular sequence generating the ideal $I$ and $
y_i\in D_{F_*}(I/I^2[1]) $ be the dual of $\bar x_i \in I/I^2[1]$.

Consider the $R_*$-algebra homomorphisms
\[
\chi_i\: (R/x_i)^*_R(R/x_i) \lra F^*_R(F)
\]
defined by $f\mapsto f\smash 1$, where $1$ denotes the identity map on $F_i'=
\smash_{j\neq i} R/x_j$ and where we identify $F$ with $R/x_i\smash F_i'$.

\begin{lem}\label{relatederivations}
Let $F=R/I$ be a regular quotient ring and let $(x_1, x_2, \ldots)$ be a regular sequence
generating the ideal $I$. For any products on $R/x_i$, $\chi_i$ restricts to an
$R_*$-homomorphism
\[
\bar\chi_i \:  \DDer^*_R(R/x_i) \lra \DDer_R^*(F).
\]
The derivations $Q_i=\bar\chi_i(Q_{x_i})$ satisfy $\psi(Q_i)= y_i$.
\end{lem}

\begin{proof}
Fix a product on $F$. To prove the first statement, it suffices to verify that
$\chi_i(\theta)\in\DDer_R^*(F)$ for $\theta\in\DDer_R^*(R/x_i)$. Choose a product $\nu$
on $R/x_i$ such that the canonical map $j\: R/x_i\to F$ is multiplicative (Proposition
\ref{natargument}). By Lemma \ref{bocksteinisderi}, $\theta$ is also a derivation with
respect to $\nu$. The diagram

\begin{equation*}
\xymatrix{ (F^\op)_*^R(R/x_i) \ar[d]^-{F^R_*(j)}\ar[r]^-{F^R_*(\theta)} &
(F^\op)_*^R(R/x_i)\ar[d]^-{F^R_*(j)}
\\
(F^\op)_*^R(F) \ar[r]^-{F^R_*(\chi_i(\theta))}
 & (F^\op)_*^R(F)  }
\end{equation*}
commutes by definition of $\chi_i(\theta)$. Since $\theta \in \DDer_R^*(R/x_i)$, Lemma
\ref{derivationlemma} implies that $F^R_*(\theta) \in \Der_{k_*}((F^\op)^R_*(R/x_i))$. We
set as usual $a_{x_i}=\phi^{F^\op}_{R/x_i}(\bar x_i) \in (F^\op)_*^R(R/x_i)$ and
$a_{j}=\phi^{F^\op}_{F}(\bar x_j) \in (F^\op)_*^R(F)$. Then $(F^\op)_*^R(R/x_i)\cong
\Lambda_{F_*}(a_{x_i})$ and $(F^\op)_*^R(F)\cong \Lambda_{F_*}(a_{1},a_{2}, \ldots)$.
Since $j$ is multiplicative and the characteristic homomorphism $\phi$ is natural,
$F^R_*(j)$ is an algebra morphism such that $F^R_*(j)(a_{x_i})=a_{i}$. Via the
isomorphism
\[
(F^\op)_*^R(F)\cong \Bigl(\bigotimes_{k\neq
i}\Lambda_{F_*}(a_{k})\Bigr)\otimes \Lambda_{F_*}(a_{i}),
\]
$F^R_*(\chi_i(\theta))$ corresponds to $1 \otimes F^R_*(\theta)$. It follows that
$F^R_*(\chi_i(\theta))$ is a derivation. By Lemma \ref{derivationlemma}, $\chi_i(\theta)$
is a derivation as well. In addition, we have $\psi(Q_i)=y_i$, by naturality of $\psi$
and by Lemma \ref{bocksteinisderi}.
\end{proof}

\begin{defn}\label{defbockstein}
Let $F=R/I$ be a regular quotient ring. The Bockstein operation
$Q_\alpha\in\DDer_R^*(F)$ associated to $\alpha\in
\Hom^*_{F_*}(I/I^2[1], F_*)$ is defined to be $\psi^{-1}(\alpha)$.
We write $Q_i$ for $Q_{y_i}$, where $(x_1,x_2, \ldots)$ is a
regular sequence generating $I$ and where $y_i$ is dual to $\bar
x_i$.
\end{defn}

\begin{rem}\label{stricklandcitation}
Strickland defines in \cite{strickland} for $F=R/I$ a regular quotient ring a
homomorphism $d\: \DDer_R^*(F)\to \Hom^*_{R^*}(I/I^2, F_*)$ and shows that $d$ is
injective. Moreover, he proves that $d$ is an isomorphism for diagonal $F$. The
homomorphism $d$ coincides with our $\psi$, as $d(Q_i)=y_i$ \cite{strickland}*{Corollary
4.19}.
\end{rem}

\begin{proof}[Proof of Proposition \ref{derivations}]
(i) We first show that $\psi$ is surjective. Choose a regular
sequence $(x_1, x_2, \ldots)$ generating $I$. Let $Q_i$ and $y_i$
be as above. By Lemma \ref{relatederivations} and by naturality of
$\psi$, we have $\psi(\pi \circ Q_i)=\pi_*\circ y_i$. Because
$\Hom^*_{F_*}(I/I^2[1], k_*)$ is generated by the elements
$\pi_*\circ y_i$, $\psi$ is surjective.

To show that $\psi$ is injective, suppose that $\theta\in\DDer_R^*(F,k)$ satisfies
$\psi(\theta)=0$. By Corollary \ref{algffop} $(\mu_k)_*\: (k^{\op})_*^R(k) \to k_*$ is
the augmentation of an exterior algebra and hence an algebra homomorphism. Therefore, the
composition
\begin{equation}\label{longcomp}
\Lambda(I/I^2[1]) \xra[\cong]{\Phi} (F^\op)_*^R(F) \xra{F^R_*(\theta)} (F^\op)_*^R(k)
\xra{(\pi\smash 1)_*} (k^\op)_*^R(k) \xra{(\mu_k)_*} k_*,
\end{equation}
where $\Phi$ is the isomorphism from Corollary \ref{algffop}, is a derivation. By
assumption, its restriction to $I/I^2[1]$ is zero. This implies that \eqref{longcomp} is
zero. By duality (see Section \ref{notations}), it follows that $\theta$ is trivial.

It remains to prove that $\psi$ is open. By definition of the topology on
$\DDer_R^*(F,k)$ and the fact that $\Psi$ is a homeomorphism (Proposition
\ref{cohomologyiso}), it suffices to show that
\[
\iota^*\: \Hom^*_{F_*}(\Lambda(I/I^2[1]), k_*) \lra
\Hom^*_{F_*}(I/I^2[1], k_*)
\]
is open. By definition of the topologies involved here, this is a consequence of the fact
that an injection of $k_*$-modules $V_*\to W_*$ induces an open map on the duals with
respect to the dual-finite topologies.

(ii) This is clear, because $\psi(\pi \circ Q_i)=\pi_*\circ y_i$
and because $Q_i=\chi_i(Q_{x_i})$ is defined independently on any
products.
\end{proof}

\begin{rem}\label{derioff}
It is a consequence of Proposition \ref{derivations} that the $\bar\chi_i$ from Lemma
\ref{relatederivations} induce an isomorphism
\[
\prod_{i\geq 1} F^*\otimes_{R^*/x_i}  \DDer_R^*(R/x_i) \cong \DDer_R^*(F).
\]
\end{rem}

We close this section by giving two properties of derivations which we will need later
on.

\begin{lem}\label{squarezero}
Any derivation $\theta\in\DDer^*_R(F)$ satisfies $\theta^2=0$.
\end{lem}

\begin{proof}
By Proposition \ref{derivations}, we may assume that $\theta=Q_i$.
By Lemma \ref{relatederivations}, we have
$\theta=\bar\chi_i(Q_{x_i})$. Hence $\theta^2$ is given by
smashing $Q_{x_i}^2$ with the identities on the other smash
factors. But $Q_{x_i}^2$ is trivial, by definition.
\end{proof}

\begin{lem}\label{effder}
For any $\theta\in\DDer_R^*(F)$, the diagram below commutes:
\[
\xymatrix{{I/I^2}[1] \ar[d]_-{\phi} \ar[r]^-{\psi(\theta)} & k_*
\ar[d]^-{(1\smash\eta)_*}
\\
k_*^R(F) \ar[r]^-{k_*^R(\theta)} & k_*^R(F).}
\]
\end{lem}

\begin{proof}
Let $\bar x\in I/I^2[1]$. Choose a product on $R/x$ such that the
canonical map $j\: R/x\to F$ is multiplicative (Proposition
\ref{natargument}). Consider the maps
\begin{equation*}
\DDer^*_R(F) \lra \DDer^*_R(R/x, F) \longleftarrow
F_*\otimes_{R_*/x} \DDer_R^*(R/x)
\end{equation*}
induced by $j$. It follows from Proposition \ref{derivations}(i)
that the second map is an isomorphism. Therefore, there is a
derivation $\theta_x\in\DDer_R^*(R/x)$ such that
\[
\xymatrix{R/x\ar[d]^-j \ar[r]^-{\theta_x} & R/x\ar[d]^-j \\ F
\ar[r]^-{\theta} & F}
\]
commutes. By naturality of $\phi$ and $\psi$, we may therefore assume that $F=R/x$. By
Proposition \ref{derivations}(i) and Lemma \ref{bocksteinisderi}, we can restrict to
$\theta=Q_x$. But then, the statement comes down to the statement in Remark
\ref{Qinduced}.
\end{proof}

\subsection{Cohomology of regular quotients}\label{cohalg}

We now determine the cohomology algebra $F^*_R(F)$ for a regular quotient $F$.

We need a notation. Let $M_*$ be a module over a graded ring
$F_*$. The dual-finite filtration on $D(M_*)=D_{F_*}(M_*)$ induces
a filtration of the exterior algebra $\Lambda(D(M_*))$. We write
$\wh\Lambda(D(M_*))$ for the completion of $\Lambda(D(M_*))$ with
respect to this filtration.

An isomorphism of the form below was constructed by Strickland for
diagonal $F$ \cite{strickland}*{Cor.\@ 4.19}. His construction
relies upon the choice of a regular sequence generating $I$. We
show that there is an isomorphism which is independent on any
choices, for any regular quotient ring $F$.

\begin{thm}\label{cohomologyalgebra}
For a regular quotient ring $F=R/I$, there is a canonical
homeomorphism of $F^*$-algebras
\[
\Theta\: \wh\Lambda(\DDer^*_R(F)) %\cong\wh\Lambda(\DDer_R^*(F))
\cong F^*_R(F).
\]
\end{thm}

\begin{rem}
Proposition \ref{derivations} and Theorem \ref{cohomologyalgebra}
imply that if $F=R/I$ is a regular quotient module, then
\[
\wh\Lambda(D(I/I^2[1])) \cong F^*_R(F).
\]
Note that on fixing a regular sequence $(x_1,x_2, \ldots)$ generating $I$, we obtain
\[
\wh\Lambda(Q_1,Q_2,\ldots) \cong F^*_R(F),
\]
where the $Q_i$ are defined according to Definition
\ref{defbockstein}.
\end{rem}

\begin{proof}
Set $V=I/I^2[1]$ and recall that $V$ is a free $F_*$-module with
basis $\bar x_1, \bar x_2, \ldots$, where $(x_1, x_2, \ldots)$ is
a regular sequence generating $I$. We define
$$
\delta \: D(V) \lra D(\Lambda(V))
$$
by $\delta(y_i)= \epsilon \circ \tfrac{\partial}{\partial \bar x_i} $ where $\epsilon$ is
the canonical augmentation of $\Lambda(V) $ and $y_i$ is dual to $\bar x_i$. We easily
check that $\delta $ lifts to a homeomorphism
$$
\Delta \: \wh\Lambda(D(V)) \lra D(\Lambda(V))
$$
with $\Delta(y_{i_1} \wedge \cdots \wedge y_{i_n})= \epsilon \circ
\tfrac{\partial}{\partial \bar x_{i_1}}\circ \cdots \circ
\tfrac{\partial}{\partial \bar x_{i_1}} $ (for the proof, consider
first the case where $V$ is finitely generated and then pass to
limits).

The homeomorphism $\Psi \: F^*_R(F) \lra D(\Lambda(V))$ from Proposition
\ref{cohomologyiso} is, in the case we are considering, just the composition of the usual
Kronecker homomorphism with the dual of the isomorphism  $\Phi \: \Lambda(V) \cong
F^R_*(F^\op)$. The Kronecker homomorphism is a homeomorphism, since $F^R_*(F^\op)$ is
$F_*$-free.

Lemma \ref{relatederivations} implies that $F^R_*(Q_i)$ is a derivation of the algebra
$F_*^R(F^\op)$. Using Remark \ref{Qinduced} and the isomorphism $F^R_*(F^\op) \cong
\Lambda (V)$, we easily check that $\Psi(Q_i)=\epsilon \circ \tfrac{\partial}{\partial
\bar x_i}$. Since $\psi(Q_i)=y_i$ (Lemma \ref{relatederivations}), we have that $\delta
\psi(Q_i)=\epsilon \circ \tfrac{\partial}{\partial \bar x_i} $. Therefore the following
diagram commutes:
\[
\xymatrix{F^*_R(F) \ar[r]^-{\Psi} & D(\Lambda(V)) \\
\DDer^*_R(F)\ar[u]^-{\subset} \ar[r]^-{\psi} &
D(V).\ar[u]^-{\delta}}
\]
As any derivation squares to 0 (Lemma \ref{squarezero}) and as $F^*_R(F)$ is complete,
the injection $\DDer^*_R(F)\hookrightarrow F^*_R(F)$ lifts to a continuous $F_*$-algebra
homomorphism
$$
\Theta \:\wh\Lambda(\DDer^*_R(F))\lra F^*_R(F).
$$
Explicitly, $\Theta$ is given by $\Theta (Q_{i_1} \wedge \cdots \wedge Q_{i_n})=Q_{i_1}
\circ \cdots \circ Q_{i_n}$. Because of $\Psi(Q_{i_1} \circ \cdots \circ
Q_{i_n})=\epsilon \circ \tfrac{\partial}{\partial \bar x_{i_1}}\circ \cdots \circ
\tfrac{\partial}{\partial \bar x_{i_1}}$, the diagram below commutes, too:
\[
\xymatrix{F^*_R(F) \ar[r]^-{\Psi} & D(\Lambda(V)) \\
\wh\Lambda(\DDer^*_R(F))\ar[u]^-{\Theta}
\ar[r]^-{\wh\Lambda(\psi)} & \wh\Lambda(D(V)).\ar[u]^-{\Delta}}
\]
Together with $\Psi$, $\Delta$ and $\wh\Lambda(\psi)$, $\Theta$ is therefore a
homeomorphism, too.
\end{proof}

\section{Examples}\label{examples}

In this section we discuss the example of the Morava $K$-theories $K(n)$.
Their 2-periodic versions $K_n$
can be treated similarly.
They are discussed in
detail in \cite{jw}.

\subsection{Definition of Morava $K$-theory}

We fix a prime number $p$. Recall that the $p$-localization
$MU_{(p)}$ of the spectrum associated to the complex cobordism
$MU$  is a commutative $\SS$-algebra (see \cite{ekmm}) satisfying:
\[
(MU_{(p)})_* \cong \Z_{(p)}[x_1,x_2, \ldots],\; |x_i|=2i.
\]

The Hopkins-Miller theorem \cite{g-h}
has as a consequence
that
for $n\geq 0$, there exists an $MU_{(p)}$-algebra $\wh
E(n)$  with
\[
\widehat{E}(n)_* \cong \lim_k
\Z_{(p)}[v_1,\ldots,v_{n-1}][v_n,v_n^{-1}]/I_n^k,
\]
where $I_n$ is the ideal generated by the regular sequence $(v_0=p, v_1, \ldots,
v_{n-1})$. Details can be found in \cite{rognes}*{Theorem 1.5} and in the unpublished
correction {\it ``A not necessarily commutative map''}, available on the author's home
page.

The $n$-th Morava $K$-theory may
be defined as the regular quotient of $\widehat{E}(n)$ by $I_n$:
\[
K(n) = \widehat{E}(n)/ I_n \cong \widehat{E}(n)/v_0
\smash_{\widehat{E}(n)} \cdots \smash_{\widehat{E}(n)}
\widehat{E}(n)/v_{n-1}.
\]
Its coefficient ring satisfies $K(n)_*\cong \F_p[v_n,v_n^{-1}].$

\subsection{The case $p$ odd}

We first consider the case where $p$ is an odd prime. According to Strickland
\cite{strickland}*{Cor.\@ 3.12}, there is a commutative $\widehat{E}(n)$-product $\mu_k$
on $\widehat{E}(n)/v_k$ for $0 \leq k \leq n-1$. Let $\mu$ be the smash ring product of
the $\mu_k$ on $K(n)$.
Since $\mu$ is commutative, we have $b_{K(n)}=0$ by Corollary
\ref{bilcomm}. Therefore if $K(n)$
is endowed with this product $\mu$, then
\[
K(n)_*^{\widehat{E}(n)}(K(n)) \cong \Lambda(I_n/I_n^2[1]) \cong
\Lambda (a_0,\ldots , a_{n-1})
\]
where $a_i=\phi(\bar v_i)$, as in Section \ref{algstructure}.

\subsection{The case $p=2$}
The case of the prime  $p=2$ is much more interesting. We use some arguments and notation
from \cite{strickland}*{Section 7} in the following.

Let $w_k\in MU_{2(2^k-1)}$ denote the bordism class of a smooth hypersurface $W_{2^k}$ of
degree $2$ in $\C P^{2^k}$ and let $J_k\subseteq (MU_{(2)})_*$ be the ideal $(w_0,
\ldots, w_{k-1})$, where $w_0=2$. The sequence of the $w_i$ is regular, and the image of
$J_k$ in $\widehat{E}(n)_*$ is the ideal $I_k=(v_0, \ldots, v_{k-1})$, for $k=0,\ldots,n$
(see \cite{strickland}). To simplify the notation, we write again $w_k$ for the image of
$w_k\in (MU_{(2)})_{*}$ in $\widehat{E}(n)_*$.

\begin{prop}\label{compofc}
There is a product $\mu_k$ on $\widehat{E}(n)/w_k$ with
$c(\mu_k)\equiv w_{k+1} \mod I_k$ for $k\geq 0$.
\end{prop}

\begin{proof}
As  $\widehat{E}(n)$ is an $MU_{(2)}$-algebra, the functor $\mathscr{E} \:
\mathscr{D}_{MU_{(2)}} \to \mathscr{D}_{\widehat{E}(n)}$ defined as $\mathscr{E}(M)=M
\wedge_{MU_{(2)}}\widehat{E}(n)$ is strictly monoidal. This can be seen as follows: For
$MU_{(2)}$-modules $M$ and $N$, $M \wedge_{MU_{(2)}}\widehat{E}(n)\cong \widehat{E}(n)
\wedge_{MU_{(2)}}M$ is an $(\widehat{E}(n),\widehat{E}(n))$-bimodule and exactly as in
\cite{ekmm}*{III. 3} there is a
natural isomorphism
\[
(M \wedge_{MU_{(2)}}\widehat{E}(n))\wedge_{\widehat{E}(n)}
(\widehat{E}(n) \wedge_{MU_{(2)}}N) \cong
\widehat{E}(n) \wedge_{MU_{(2)}}(M \wedge_{MU_{(2)}}N)
\]
of $(\widehat{E}(n),\widehat{E}(n))$-bimodules.
As a consequence, the functor
$\mathscr{E}$
maps $MU_{(2)}$-rings to $\widehat{E}(n)$-rings. Strickland
constructs a $MU_{(2)}$-product $\tilde\mu_k$ on $MU_{(2)}/w_k$
with $c(\tilde\mu_k)\equiv w_{k+1} \mod J_k$ for $k\geq 0$
\cite{strickland}*{Section 7}. Via the functor $\mathscr{E}$,
$\tilde\mu_k$ induces an $\widehat{E}(n)$-product $\mu_k$ on
$\widehat{E}(n)/w_k$. By definition of the obstruction $c$, we
check that $c(\mu_k)\equiv w_{k+1} \mod I_k$.
\end{proof}

We endow $K(n)$ with the diagonal product $\mu$, defined as the smash ring product of the
$\mu_k$. As $v_n  \equiv w_n \mod I_n$, Propositions \ref{bildiag} and \ref{compofc}
imply that $b_{K(n)}= v_n \cdot y_{n-1} \otimes y_{n-1}$, where $y_{n-1}\in
D_{\widehat{E}(n)_*}(I_n/I_n^2[1])$ is dual to $\bar v_{n-1} \in I_n/I_n^2[1]$.
Therefore, $\mu$ is not commutative, see Corollary \ref{bilcomm}.

The opposite product $\mu^\op$ is the smash ring product of the $\mu_k^\op$. It follows
from \cite{strickland}*{Prop.\@ 3.1 and Lemma\@ 3.11} that $c(\mu_k^\op)\equiv w_{k+1}
\mod I_k$, as $2 \in I_k$.

Let $1\leq k \leq n$. For dimensional reasons, we have
$(\widehat{E}(n)/w_{k-1})_{2|w_{k-1}|+2} = \{0,v_k \}$. Therefore, Proposition
\ref{basicproducts} implies that
\[
\mu_{k-1}^\op=\mu_{k-1} \circ (1+v_k \cdot Q_{w_{k-1}} \wedge Q_{w_{k-1}}).
\]
The elements $\bar w_{k-1}, \bar v_{k-1}\in I_k/I_k^2[1]$ coincide, hence their duals are
the same and so Proposition \ref{derivations} implies that $Q_{w_{k-1}} = Q_{v_{k-1}} \in
\DDer_{\widehat{E}(n)}^*(\widehat{E}(n)/w_{k-1})$. As a consequence, we recover the well
known formula:
\[
\mu^\op=\mu \circ (1+v_n \cdot Q_{n-1} \wedge Q_{n-1}),
\]
where $Q_{n-1}$ is defined as in Definition \ref{defbockstein}.
Observe that $b_{K(n)}=b_{K(n)^\op}$ although $\mu \neq \mu^\op$.
We now compute
\[
K(n)_*^{\widehat{E}(n)}(K(n)) \cong \Lambda(a_0,
\ldots,a_{n-2})\otimes \T(a_{n-1})/( a_{n-1}^2-v_n \cdot 1),
\]
where $K(n)$
is endowed with the product $\mu$
described above and
the $a_i$ are defined as in the case where $p$ is odd.

\appendix

\section{An algebraic fact}\label{commalg}

\begin{prop}
Suppose that $I_1, I_2, \ldots\subset R_*$ are ideals which
satisfy
\[
(I_1+ \cdots + I_{k-1}) \cdot I_k = (I_1 + \cdots + I_{k-1})\cap
I_k
\]
for all $k>1$. Let $I=I_1+I_2+\cdots$. Then there is a canonical isomorphism of
$R_*/I$-modules
\begin{equation*}%\label{decompdiag}
I/I^2\cong \bigoplus_{i\geq 1}  R/I_* \otimes_{R_*} I_i/I_i^2.
\end{equation*}
\end{prop}

\begin{proof}
We prove the statement only for $I = I_1 + I_2$. The argument
needed for the inductive step is similar and therefore left to the
reader. For infinitely many $I_i$, the statement follows by
passing to colimits.

We begin by showing that $I^2 = (I_1+I_2^2)\cap (I_1^2 + I_2)$. The inclusion $\subseteq$
is trivial. To show $\supseteq$, suppose that $\alpha\in(I_1+I_2^2)\cap (I_1^2 + I_2)$.
Write $\alpha$ as $\alpha = x + w = v + y$, where $x\in I_1$, $w\in I_2^2$, $y\in I_2$
and $v\in I_1^2$. It follows that $x-v=y-w\in I_1\cap I_2$. By hypothesis, we have
$I_1\cap I_2 = I_1\cdot I_2$, and therefore $\alpha= (x-v) + v + w \in I_1\cdot I_2 +
I_1^2 + I_2^2 =I^2$.

It follows that the canonical homomorphism
\begin{equation}\label{aniso}
I/I^2 \lra I/(I_1^2 + I_2) \oplus I/(I_1+I_2^2)
\end{equation}
is an isomorphism. Moreover, the canonical map
\begin{equation}\label{anotheriso}
I_1/(I_1\cap I_2 + I_1^2) \lra I/(I_1^2 +I_2)
\end{equation}
and its symmetric analogue are easily seen to be isomorphisms. Finally, there is a
natural isomorphism
\begin{equation}\label{lastiso}
R_*/I \otimes_{R_*} I_1/I_1^2 \cong I_1/(I_1\cap I_2 + I_1^2),
\end{equation}
given by the following composition:
\begin{align*}
R_*/I \otimes_{R_*} I_1/I_1^2 & \cong R_*/I_2\otimes_{R_*}
I_1/I_1^2 \cong (I_1/I_1^2)/\bigl(I_2\cdot(I_1/I_1^2)\bigr)
\\
& \cong I/(I_1 \cdot I_2+ I_1^2) \cong I_1/(I_1\cap I_2 + I_1^2).
\end{align*}
Combining \eqref{aniso}, \eqref{anotheriso} and \eqref{lastiso}
implies the result.
\end{proof}

\end{document}